\documentclass{amsart}

\usepackage{geometry,graphicx,amssymb,amsmath,amsbsy,eucal,amsfonts,mathrsfs,amscd,bm,tcolorbox, xcolor,caption,float,tikz,tikz-3dplot,cancel}
\usepackage[font=small]{subcaption}
\usetikzlibrary{patterns}
\usepackage[textsize=tiny,color=blue!50!white]{todonotes}
\usepackage[normalem]{ulem} % for strikethrough (\sout) in collaborative writing

\usepackage[breaklinks,bookmarks=false]{hyperref}

\hypersetup{colorlinks, linkcolor=blue, citecolor=blue,
urlcolor=blue, plainpages=false, pdfwindowui=false,
pdfstartview={FitH}, pdftitle={},}

\numberwithin{equation}{section}

\allowdisplaybreaks[3]

\newtheorem{theorem}{Theorem}[section]
\newtheorem{lemma}[theorem]{Lemma}

\newtheorem{proposition}[theorem]{Proposition}
\theoremstyle{definition}

\newtheorem{remark}[theorem]{Remark}

\newcommand{\rev}[1]{{\color{black}#1}}

\usepackage{mydef}

\begin{document}

\title[Local $L^2$-bounded commuting projections using discrete problems]{Local $L^2$-bounded commuting projections using discrete local problems on Alfeld splits}

\def\AEr{Alexandre Ern}
\def\AEa{CERMICS, ENPC, Institut Polytechnique de Paris, 77455 Marne-la-Vall\'ee, France \& Inria Paris, 48 rue Barrault, 75647 Paris, France}

\def\JG{Johnny Guzm\'an}
\def\JGa{Division of Applied Mathematics,
Brown University,
Box F,
182 George Street,
Providence, RI 02912, USA}

\def\PP{Pratyush Potu}
\def\PPa{Division of Applied Mathematics,
Brown University,
Box F,
182 George Street,
Providence, RI 02912, USA}

\def\MV{Martin Vohral\'ik}
\def\MVa{Inria Paris, 48 rue Barrault, 75647 Paris, France \& CERMICS, ENPC, Institut Polytechnique de Paris, 77455 Marne-la-Vall\'ee, France}

\author{\AEr}
\address{\AEa}
\email{\href{mailto:alexandre.ern@enpc.fr}{alexandre.ern@enpc.fr}}

\author{\JG}
\address{\JGa}
\email{\href{mailto:johnny_guzman@brown.edu}{johnny\_guzman@brown.edu}}

\author{\PP}
\address{\PPa}
\email{\href{mailto:pratyush_potu@brown.edu}{pratyush\_potu@brown.edu}}

\author{\MV}
\address{\MVa}
\email{\href{mailto:martin.vohralik@inria.fr}{martin.vohralik@inria.fr}}

%\thanks{}


\begin{abstract}
We construct projections onto the classical finite element spaces based on Lagrange, N\'ed\'elec, Raviart--Thomas, and discontinuous elements on shape-regular simplicial meshes. 
Our projections are defined locally, are bounded in the $L^2$-norm, and commute with the corresponding differential operators.
\rev{Such projections are a fundamental tool in finite element stability and error analysis. Moreover, to the best of our knowledge, the present construction is the first in the literature where local $L^2$-stability is fully established.}
The cornerstone of our construction are local weight functions which are piecewise polynomials built using the Alfeld split of local patches from
the original simplicial mesh. 
This way, the $L^2$-stability of the projections is established by invoking discrete
Poincar\'e inequalities on these local stars, for which we 
provide \rev{here an original, constructive proof}. \rev{Another important novelty is that
we extend the construction of the projections so as
to preserve homogeneous boundary conditions on a subset of the domain boundary}. 
Altogether, the material is presented using the language of vector calculus, and links to the formalism of finite element exterior calculus are provided.
\end{abstract}

\maketitle

% \setcounter{tocdepth}{3}
% \tableofcontents

%%%%%%%%%%%%%%%%%%%%%%%%%%%%%%%%%%%%%
%%%%%%%%%%%%%%%%%%%%%%%%%%%%%%%%%%%%%

\section{Introduction}

Bounded commuting projections are a key tool in the finite element exterior calculus (FEEC). 
\rev{Indeed, they are needed in the finite element stability and error analysis 
of primal and mixed formulations of PDEs involving the curl and/or the divergence operators.
Such projections conceptually differ from 
the so-called canonical interpolation operators of Raviart and Thomas~\cite{Ra_Tho_MFE_77}, Fortin~\cite{Fort_MFEs_77}, and N\'ed\'elec~\cite{nedelec1980mixed}, for which the required regularity is much stronger than the one provided by the PDE at hand.}
Some of the first bounded commuting projections were derived by Sch\"oberl~\cite{schoberlmultilevel} and Christiansen and Winther~\cite{Christ_Wint_sm_proj_08} \rev{and were globally defined}. 
Later, Falk and Winther~\cite{Falk_Winth_loc_coch_14,falk2015double} constructed commuting projections that were in addition local (i.e., to define the projection on a simplex, only information on a local neighborhood of the simplex is needed); however, these projections are bounded in the corresponding graph norm and not in $L^2$.
Inspired by this work, Arnold and Guzm\'an~\cite{Arn_Guz_loc_stab_L2_com_proj_21} constructed commuting projections that are local and whose operator norm is locally bounded in $L^2$, under a conjecture discussed in more detail below. It is also worth mentioning that a different approach to construct local commuting projections has been carried out by Ern et al.~\cite{Ern_Gud_Sme_Voh_loc_glob_div_22} and Chaumont-Frelet and Vohral\'ik~\cite{Chaum_Voh_H_curl_proj_24}\rev{, where local $L^2$-boundedness holds up to ``data oscillation'' on the curl or the divergence. Altogether, we see that a theoretical gap remains in the literature to devise local commuting projections for which local $L^2$-stability is fully established. The main contribution of the present work is to close this gap. $L^2$-stability is, in particular, important to derive improved error estimates (\cite[Theorem~3.11]{arnold2010finite}) and convergence for eigenvalue problems (\cite[Theorem~3.19]{arnold2010finite}).}

The underlying functional setting for the commuting projections under consideration is
the well-known de Rham complex
\begin{alignat}{4}\label{complex}
&\mathbb{R}
\stackrel{\subset}{\xrightarrow{\hspace*{0.5cm}}}\
 V^0
&&\stackrel{\grad}{\xrightarrow{\hspace*{0.5cm}}}\
 \bV^1
&&\stackrel{\curl}{\xrightarrow{\hspace*{0.5cm}}}\
 \bV^2
&&\stackrel{\dive}{\xrightarrow{\hspace*{0.5cm}}}\
V^3
\stackrel{}{\xrightarrow{\hspace*{0.5cm}}}
0,
\end{alignat}
where the relevant graph spaces are
\begin{equation} \label{eq:graph_spaces}
V^0:= H(\grad,\Omega) =H^1(\Omega), \quad \bV^1:= \bH(\curl, \Omega), \quad \bV^2:= \bH(\dive, \Omega), \quad V^3:= L^2(\Omega). 
\end{equation}
The key idea in~\cite{Arn_Guz_loc_stab_L2_com_proj_21} is to build weight functions that have commuting properties with the differential operators in~\eqref{complex} (also called exterior derivatives in the FEEC context) and, at the same time, have local bounds in $L^2$. 
To do this, the authors in~\cite{Arn_Guz_loc_stab_L2_com_proj_21} solved local problems (on extended stars of simplices) and then invoked the work of Costabel and McIntosh~\cite{Cost_McInt_Bog_Poinc_10}, who proved the existence of a bounded right-inverse of the exterior derivative on Lipschitz domains. 
In particular, bounds on the right inverse in the $H^1$-norm can be shown if the domain is star shaped with respect to a ball of similar size. 
However, given an arbitrary simplicial mesh, it might be that the extended stars are not star shaped with respect to a ball. Nonetheless, \cite{Arn_Guz_loc_stab_L2_com_proj_21} conjectured such bounds to establish the expected bounds on the weight functions and consequently to bound the operator norm of the commuting projection in $L^2$.

The main contribution of this paper is to close the theoretical gap left in~\cite{Arn_Guz_loc_stab_L2_com_proj_21}. 
We follow~\cite{Arn_Guz_loc_stab_L2_com_proj_21} but introduce two important novel ideas. 
First, we \rev{\emph{redefine} the construction of the weight functions to be fully discrete}. We achieve this by solving local problems on piecewise polynomial spaces on a certain refinement of the actual mesh known in the literature as the Alfeld split~\cite{Alfel84}. \rev{Working on the Alfeld split is not a necessity, but a simple and convenient way to have at hand suitable bubble functions that on the one hand vanish on all the mesh faces and on the other hand are piecewise (affine) polynomials. This latter property is, in turn, crucial to invoke discrete Poincar\'e inequalities (as opposed to \rev{continuous right-inverses in $H^1$}) to establish the $L^2$-stability of the weight functions used to define the projections. Using piecewise \emph{affine} bubble functions allows us to consider the \emph{same} polynomial degree for the discrete Poincar\'e inequalities on the Alfeld split.} 
Second, we \rev{\emph{provide constructive proofs}} of the discrete Poincar\'e inequalities. % that are needed to establish the $L^2$-bounds on the operator norms. 
This latter point is in itself a result of independent interest. We include here a compact self-contained presentation \rev{with an original proof, and refer the reader to~\cite{PGEV25} for a broader discussion on discrete Poincar\'e inequalities and a review on their proofs}. 
Overall, our work leads to local commuting projections with fully provable local $L^2$-bounds on their operator norm. To our knowledge, this is the first time that projections with all such properties are established in the literature. 

\rev{The last salient contribution of our work is to show how to extend the construction so as to prescribe \emph{homogeneous boundary conditions} on a subset of the boundary. This additional property, which is not addressed in \cite{Arn_Guz_loc_stab_L2_com_proj_21}, is crucial in applications where Dirichlet or mixed Dirichlet/Neumann boundary conditions are enforced. In fact, the projections of \cite[Section 6]{Christ_Wint_sm_proj_08} (see also \cite[Chapter~23]{ErnGuermondbook}) preserve homogeneous boundary conditions on the entire boundary, while similar globally defined $L^2$-bounded commuting projections which preserve homogeneous boundary conditions on a subset of the boundary have been devised by Gopalakrishnan and Qiu \cite[Section 3.2]{Gop_Qiu_Lip_exp_12} and Licht \cite{Licht_proj_bc_19}. Homogeneous boundary conditions on a subset of the boundary are also preserved by the projectors of~\cite{Ern_Gud_Sme_Voh_loc_glob_div_22, Chaum_Voh_H_curl_proj_24}, which are locally $L^2$-stable up to data oscillation, and by those of Hiptmair and PEchstein \cite[Section 3.2.6]{Hipt_Pech_disc_1fo_19}, which are stable in the graph norm.}

As a final remark on the literature, we notice that the weight functions in~\cite{Falk_Winth_loc_coch_14,falk2015double} are also discrete, but two correction steps are required in their construction (as opposed to one here and in~\cite{Arn_Guz_loc_stab_L2_com_proj_21}). 
Moreover, we achieve local $L^2$-stability of the projections, whereas the projections in~\cite{Falk_Winth_loc_coch_14,falk2015double} only enjoy stability in the corresponding graph norms. 

The paper is organized as follows. 
In Section~\ref{prelim}, we introduce the discrete setting. 
In Section~\ref{sec:weight}, we present the requested properties of the novel weight functions. The main result of this section is Theorem~\ref{th:weights}.
The actual construction of the weight functions and the proof of their key properties 
are collected in Section~\ref{sec:proof}. \rev{Discrete Poincar\'e inequalities are needed here, and we present a constructive proof in the appendix.} 
In Section~\ref{sec:projections}, we define the local $L^2$-bounded commuting projections.
Finally, the extensions needed
so that the projections additionally satisfy homogeneous conditions \rev{on a subset of} the boundary are
discussed in Section~\ref{sec:boundary}.
Although it is possible to prove all the results in the paper in any dimension 
using the language of 
FEEC spaces~\cite{arnold2006finite, arnold2010finite, arnold2018finite},
we have chosen here to focus on three space dimensions and to work with vector notation, as in~\eqref{complex}--\eqref{eq:graph_spaces}. We hope that
this presentation will make the material more accessible to a wider audience and allow readers to appreciate more thoroughly the results here and in~\cite{Arn_Guz_loc_stab_L2_com_proj_21}.

\section{Discrete setting}\label{prelim}

In this section, we present the discrete setting, namely the discrete objects 
associated with the mesh and the piecewise polynomial spaces
based on Lagrange, N\'ed\'elec, Raviart--Thomas, and discontinuous finite elements.
We also define %traces on the discrete geometric objects, as well as
the canonical degrees of freedom together with the Whitney forms associated with 
the lowest-order version of the above spaces.

\subsection{Simplicial mesh}
Let $\Omega$ be a Lipschitz, polyhedral, open, bounded, connected set in $\RRR^3$. 
Let $\Th$ be a simplicial triangulation of $\Omega$.
For all $l\in\{0{:}3\}$, the $l$-simplices in $\Th$ are the mesh vertices for $l=0$,
the mesh edges for $l=1$, the mesh faces for $l=2$, and the mesh tetrahedra for $l=3$.
Notice that $l$-simplices are, by definition, closed sets.
We enumerate the vertices of $\Th$ and denote the set of vertices by $\Vh:=\{ x_0, \ldots, x_{N}\}$.
All the $l$-simplices are oriented by taking their vertices in
increasing enumeration order \cite[Sec.~10.3]{ErnGuermondbook}.
We denote the collection of (oriented) $l$-simplices as
\begin{equation} \label{eq:def_sigma}
\Delta_h^l := \{ \sigma = [x_{i_0},\ldots,x_{i_l}] \,:\, 0\le i_0<\ldots<i_l\le N\},
\end{equation}
where the brackets denote the convex hull of a set of points.
A more explicit notation is
\begin{equation}
\Vh := \Delta_h^0, \qquad \Eh := \Delta_h^1, \qquad \Fh := \Delta_h^2, \qquad \Th := \Delta_h^3,
\end{equation}
where $\Eh$ is the set of (oriented) mesh edges, $\Fh$ is the set of (oriented) mesh faces, and $\Th$ the set of (oriented) mesh cells
(tetrahedra). 
The shape-regularity parameter of the mesh $\Th$ is defined as 
\begin{equation} \label{eq_rho}
\rho_{\Th} \eq \max_{\tau \in \Th} h_\tau / \iota_\tau,
\end{equation}
where $h_\tau$ is the diameter of $\tau$ and $\iota_\tau$ the diameter of the largest ball inscribed in $\tau$.

The set $\Delta_h := \bigcup_{l\in\{0{:}3\}} \Delta_h^l$ is the collection of all the (oriented) geometric objects in the mesh. 
For every edge $e:=[x_{i_0}, x_{i_1}] \in \Eh$, we let $\bt_e$ be the unit tangent vector to $e$ pointing from $x_{i_0} $ to $x_{i_1}$. 
For every face $f:= [x_{i_0}, x_{i_1}, x_{i_2}] \in \Fh$, we let $\bn_f$ be the unit normal vector to $f$ such that $\bn_f:=\bt_{e_1}{\times}\bt_{e_2}$
with $e_1:=[x_{i_0}, x_{i_1}]$ and $e_2:=[x_{i_0}, x_{i_2}]$.
It is convenient to define the subsets 
\begin{subequations} \begin{alignat}{2}
\Ve&:=\{v\in\Vh\,:\, v\in e\}, &&\quad \forall e\in \Eh, \\
\Vf&:=\{v\in\Vh\,:\, v\in f\}, &&\quad \forall f\in \Fh, \\
\Vt&:=\{v\in\Vh\,:\, v\in \tau\}, &&\quad \forall \tau\in \Th, \\
\Ef&:=\{e\in\Eh\,:\, e\subset f\}, &&\quad \forall f\in \Fh, \\
\Et&:=\{e\in\Eh\,:\, e\subset \tau\}, &&\quad \forall \tau\in \Th, \\
\Ft&:=\{f\in\Fh\,:\, f\subset \tau\}, &&\quad \forall \tau\in \Th.
\end{alignat} \end{subequations}

\begin{figure}[htb]
  \centering
  \begin{subfigure}[b]{0.4\textwidth}
     \centering
     \tdplotsetmaincoords{80}{165}
     \begin{tikzpicture}[tdplot_main_coords, scale = 1]
     
     % Define coordinates for the vertices of the outer tetrahedron
     \coordinate (A) at (0,0,0);
     \coordinate (B) at (4,0,0);
     \coordinate (C) at (2,3.5,0);
     \coordinate (D) at (2,1.2,4);
     
     % Draw the outer tetrahedron
     \draw[fill=blue!0, opacity=0.4] (A) -- (B) -- (C) -- cycle;
     \draw[fill=blue!0, opacity=0.4] (A) -- (B) -- (D) -- cycle;
     \draw[fill=blue!0, opacity=0.4] (A) -- (C) -- (D) -- cycle;
     \draw[fill=blue!0, opacity=0.4] (B) -- (C) -- (D) -- cycle;
     
     % Define coordinates for the vertices of the inner tetrahedron
     \coordinate (E) at (1.5, 1, 1);
     \coordinate (F) at (2.5, 1, 1);
     \coordinate (G) at (2, 2.5, 1);
     \coordinate (H) at (2, 1.5, 2.5);
     
     % Coloring upper tetrahedron
     \draw[fill=blue!20, opacity=0.4] (E) -- (F) -- (H) -- cycle;
     \draw[fill=blue!20, opacity=0.4] (E) -- (G) -- (H) -- cycle;
     \draw[fill=blue!20, opacity=0.4] (F) -- (G) -- (H) -- cycle;

     % Coloring bottom tetrahedron
     \draw[fill=blue!20, opacity=0.4] (E) -- (F) -- (C) -- cycle;
     \draw[fill=blue!20, opacity=0.4] (E) -- (G) -- (C) -- cycle;
     \draw[fill=blue!20, opacity=0.4] (F) -- (G) -- (C) -- cycle;

     %Coloring Face
     \draw[fill=red!60, opacity=0.8] (E) -- (F) -- (G) -- cycle;

     % Draw dashed lines connecting inner vertices to outer vertices
     \draw[dashed, thick] (F) -- (B);
     \draw[dashed, thick] (G) -- (B);
     \draw[dashed, thick] (H) -- (B);
     \draw[dashed, thick] (G) -- (C);
     \draw[dashed, thick] (F) -- (C);
     \draw[dashed, thick] (E) -- (C);
     \draw[dashed, thick] (G) -- (A);
     \draw[dashed, thick] (E) -- (A);
     \draw[dashed, thick] (H) -- (A);
     \draw[dashed, thick] (H) -- (D);
     \draw[dashed, thick] (E) -- (D);
     \draw[dashed, thick] (F) -- (D);
     \draw[dashed, thick] (F) -- (H);
     \draw[dashed, thick] (E) -- (H);
     \draw[dashed, thick] (G) -- (H);
     
     % Draw nodes for vertices
     \foreach \point in {A,B,C,D,E,F,G,H}
         \fill[black] (\point) circle (1.5pt);

     \end{tikzpicture}
     \caption{$\st(\sigma)$ where $\sigma\in\Fh$}
  \end{subfigure}
  \hspace{2mm}
  \begin{subfigure}[b]{0.4\textwidth}
     \centering

     %Set viewing angle {Altitude}{Azimuth}
     \tdplotsetmaincoords{80}{165}
     \begin{tikzpicture}[tdplot_main_coords, scale = 1]
     
     % Define coordinates for the vertices of the outer tetrahedron
     \coordinate (A) at (0,0,0);
     \coordinate (B) at (4,0,0);
     \coordinate (C) at (2,3.5,0);
     \coordinate (D) at (2,1.2,4);
     
     % Draw the outer tetrahedron
     \draw[fill=blue!20, opacity=0.4] (A) -- (B) -- (C) -- cycle;
     \draw[fill=blue!20, opacity=0.4] (A) -- (B) -- (D) -- cycle;
     \draw[fill=blue!20, opacity=0.4] (A) -- (C) -- (D) -- cycle;
     \draw[fill=blue!20, opacity=0.4] (B) -- (C) -- (D) -- cycle;
     
     % Define coordinates for the vertices of the inner tetrahedron
     \coordinate (E) at (1.5, 1, 1);
     \coordinate (F) at (2.5, 1, 1);
     \coordinate (G) at (2, 2.5, 1);
     \coordinate (H) at (2, 1.5, 2.5);
     
     % Draw the inner tetrahedron
     \draw[fill=blue!20, opacity=0.4] (E) -- (F) -- (H) -- cycle;
     \draw[fill=blue!20, opacity=0.4] (E) -- (G) -- (H) -- cycle;
     \draw[fill=blue!20, opacity=0.4] (F) -- (G) -- (H) -- cycle;
     \draw[fill=red!60, opacity=0.8] (E) -- (F) -- (G) -- cycle;
     
     % Draw dashed lines connecting inner vertices to outer vertices
     \draw[dashed, thick] (F) -- (B);
     \draw[dashed, thick] (G) -- (B);
     \draw[dashed, thick] (H) -- (B);
     \draw[dashed, thick] (G) -- (C);
     \draw[dashed, thick] (F) -- (C);
     \draw[dashed, thick] (E) -- (C);
     \draw[dashed, thick] (G) -- (A);
     \draw[dashed, thick] (E) -- (A);
     \draw[dashed, thick] (H) -- (A);
     \draw[dashed, thick] (H) -- (D);
     \draw[dashed, thick] (E) -- (D);
     \draw[dashed, thick] (F) -- (D);
     \draw[dashed, thick] (F) -- (H);
     \draw[dashed, thick] (E) -- (H);
     \draw[dashed, thick] (G) -- (H);
     
     % Draw nodes for vertices
     \foreach \point in {A,B,C,D,E,F,G,H}
         \fill[black] (\point) circle (1.5pt);
     
     \end{tikzpicture}
     \caption{$\es(\sigma)$ where $\sigma\in\Fh$}
  \end{subfigure}
     \caption{Example of $\st(\sigma)$ and $\es(\sigma)$ for $\sigma\in\Fh$. Here, $\sigma$ is highlighted in red, and the star (resp. extended star) of $\sigma$ is shaded in blue.}
     \label{fig:extended_stars}
\end{figure}
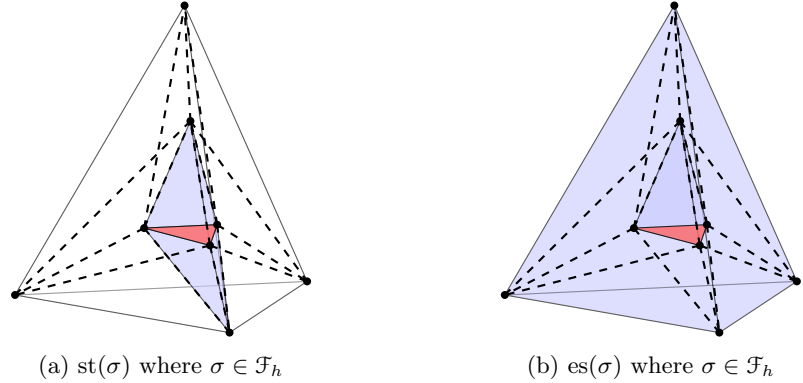

For all $l\in\{0{:}3\}$ and all $\sigma \in \Delta_h^l$, we define the \emph{star} and the \emph{extended star} of $\sigma$ as
\begin{equation*}
\st(\sigma) := \inte \bigcup_{\substack{\tau \in  \Th \\\sigma \subset \tau}} \tau, \qquad
\es(\sig) := \inte \bigcup_{\substack{\tau \in  \Th \\ \sigma \cap  \tau \neq \emptyset}} \tau.
\end{equation*}
Notice that, by definition, $\st(\sigma)$ and $\es(\sigma)$ are open subsets of $\R^3$;
we denote their closure as $\clos(\st(\sigma))$ and $\clos(\es(\sigma))$.
Equivalently, $\es(\sigma)$ is the union of $\st(v)$ for all the vertices $v\in \sigma$. In particular, if $\sigma$ is a vertex, then $\st(\sigma)=\es(\sigma)$, and if $\sigma$ is a tetrahedron, $\st(\sigma)=\sigma$. Some illustrations are shown in Figure~\ref{fig:extended_stars}.
We define $h_{\sigma}:= \diam(\sigma)$ if $l \ge 1$ and $h_{\sigma}:=\diam(\st(\sigma))$ if $l=0$. 
As in~\cite{Falk_Winth_loc_coch_14} and~\cite{Arn_Guz_loc_stab_L2_com_proj_21}, we assume that $\clos(\es(\sig))$ is contractible \cite[Page~108]{Munkres_1984} for all $\sig$ in $\Delta_h$, as is usually the case. Then, the sequences considered below in~\eqref{localcomplex} and~\eqref{localcomplexm} are exact. Alternatively, this exactness also follows if $\es(\sig)$ is simply connected with $\partial \es(\sig)$ connected \cite[Corollary~2.4, Theorem~2.9, Remark~3.10]{Gir_Rav_NS_86}. Remark~\ref{rem:non-contract} brings some further insight.

\begin{figure}[htb]
  \centering
  \begin{subfigure}[b]{0.4\textwidth}
  \begin{center}
    \tdplotsetmaincoords{80}{165}
  \begin{tikzpicture}[tdplot_main_coords, scale = 1]
  
  % Define coordinates for the vertices of the outer tetrahedron
  \coordinate (A) at (0,0,0);
  \coordinate (B) at (4,0,0);
  \coordinate (C) at (2,3.5,0);
  \coordinate (D) at (2,1.2,4);
  
  % Draw the outer tetrahedron
  \draw[fill=blue!20, opacity=0.4] (A) -- (B) -- (C) -- cycle;
  \draw[fill=blue!20, opacity=0.4] (A) -- (B) -- (D) -- cycle;
  \draw[fill=blue!20, opacity=0.4] (A) -- (C) -- (D) -- cycle;
  \draw[fill=blue!20, opacity=0.4] (B) -- (C) -- (D) -- cycle;
  
  % Define coordinates for the vertices of the inner tetrahedron
  \coordinate (E) at (1.5, 1, 1);
  \coordinate (F) at (2.5, 1, 1);
  \coordinate (G) at (2, 2.5, 1);
  \coordinate (H) at (2, 1.5, 2.5);
  
  % Draw dashed lines for inner edges
  \draw[dashed, thick] (F) -- (B);
  \draw[dashed, thick] (G) -- (B);
  \draw[dashed, thick] (H) -- (B);
  \draw[dashed, thick] (G) -- (C);
  \draw[dashed, thick] (F) -- (C);
  \draw[dashed, thick] (E) -- (C);
  \draw[dashed, thick] (G) -- (A);
  \draw[dashed, thick] (E) -- (A);
  \draw[dashed, thick] (H) -- (A);
  \draw[dashed, thick] (H) -- (D);
  \draw[dashed, thick] (E) -- (D);
  \draw[dashed, thick] (F) -- (D);
  \draw[dashed, thick] (G) -- (H);
  \draw[dashed, thick] (F) -- (H);
  \draw[dashed, thick] (E) -- (H);

  \draw[fill= gray!100, opacity=0.7] (F) -- (G) -- (H) -- cycle;
  \draw[fill= gray!100, opacity=0.7] (E) -- (G) -- (H) -- cycle;
  \draw[pattern=north west lines, pattern color=red] (F) -- (G) -- (E) -- cycle;

  % Draw nodes for vertices
  \foreach \point in {A,B,C,D,E,F,G,H}
      \fill[black] (\point) circle (1.5pt);
  
  \end{tikzpicture}
  \end{center}
  \caption{3d non-contractible $\clos(\es(\sigma))$}
\end{subfigure}
\hspace{2mm}
\begin{subfigure}[b]{0.4\textwidth}
  \begin{center}
  \begin{tikzpicture}[scale = 1]

  % Define coordinates for the vertices of the outer domain
  \coordinate (A) at (0,2);
  \coordinate (B) at (2,0);
  \coordinate (C) at (0,-1);
  \coordinate (D) at (-2,0);

  % Define coordinates for interior vertices
  \coordinate (E) at (0, 1);
  \coordinate (F) at (1, 0);
  \coordinate (G) at (-1, 0);

  % Draw the triangles in the extended star
  \draw[fill=blue!30, opacity=0.4] (A) -- (B) -- (F) -- cycle;
  \draw[fill=blue!30, opacity=0.4] (B) -- (C) -- (F) -- cycle;
  \draw[fill=blue!30, opacity=0.4] (C) -- (F) -- (G) -- cycle;
  \draw[fill=blue!30, opacity=0.4] (C) -- (D) -- (G) -- cycle;
  \draw[fill=blue!30, opacity=0.4] (A) -- (D) -- (G) -- cycle;
  \draw[fill=blue!30, opacity=0.4] (A) -- (G) -- (E) -- cycle;
  \draw[fill=blue!30, opacity=0.4] (A) -- (E) -- (F) -- cycle;

  \draw[fill= gray!100, opacity=1.0] (E) -- (F) -- (G) -- cycle;

  \draw[thick, red] (F) -- (G);

  % Draw nodes for vertices
  \foreach \point in {A,B,C,D,E,F,G}
      \fill[black] (\point) circle (1.5pt);

  \end{tikzpicture}
  \end{center}
  \caption{2d non-contractible $\clos(\es(\sigma))$}
  \end{subfigure}
s  \caption{Example showing an extended star $\es(\sigma)$ (shaded in blue) of a face $\sigma\in\Fh$ (hatched in red) where $\clos(\es(\sig))$ is non-contractible. Here, the gray region represents a hole in the domain $\Omega$. A similar 2D example (where $\sigma$ highlighted in red is now an edge) is also included for clarity.}
  \label{noncontractibleextendedstar}
\end{figure}
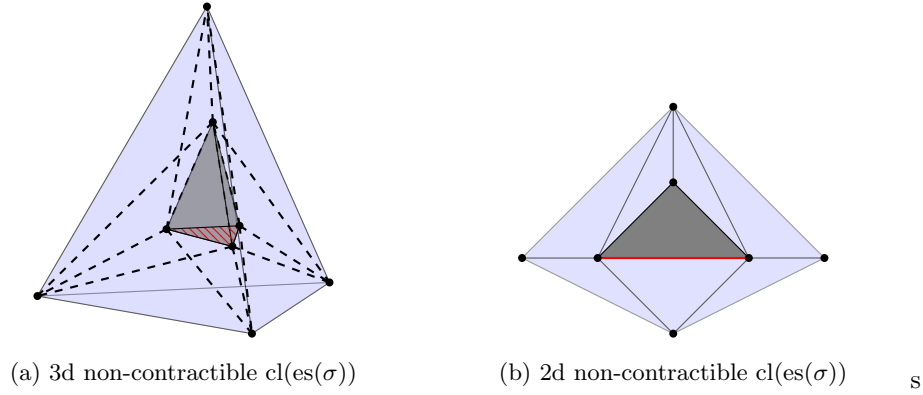

\begin{figure}[htb]
\begin{subfigure}[b]{0.45\textwidth}
\centerline{\includegraphics[width=0.7\textwidth]{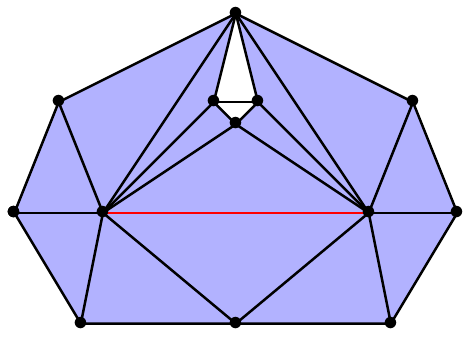}}
\caption{$\clos(\es(\sig))$ non-contractible, $\es(\sig)$ contractible}
\end{subfigure}
\begin{subfigure}[b]{0.45\textwidth}
\centerline{\includegraphics[width=0.7\textwidth]{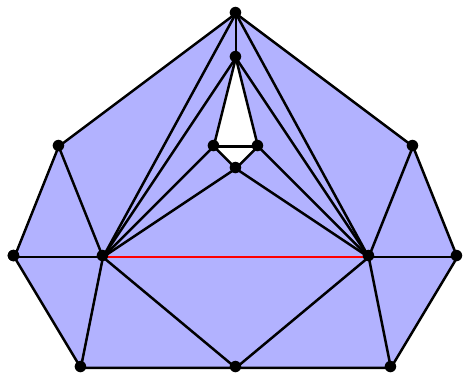}}
\caption{both $\clos(\es(\sig))$ and $\es(\sig)$ non-contractible}
\end{subfigure}
\caption{Two-dimensional illustration of a non-contractible closure of an extended star without a hole in the domain $\Omega$ (here, the white region is part of the domain $\Omega$). Subfigure (a) additionally illustrates the difference between considering the closure or not of the extended star. Extended star $\es(\sigma)$ (shaded in blue) of an edge $\sig$ (highlighted in red).}
\label{fig_ext_star_not_contractible}
\end{figure}

\begin{remark}[Non-contractible extended stars] \label{rem:non-contract}
Figure~\ref{noncontractibleextendedstar} shows examples of extended stars $\es(\sigma)$ of a face $\sigma\in\Fh$ where $\clos(\es(\sig))$ is non-contractible. Here, $\clos(\es(\sig))$ is non-contractible because of the hole in the domain $\Omega$. 
Figure~\ref{fig_ext_star_not_contractible} shows that such situations can appear without the presence of a hole in the domain $\Omega$ since the white region is part of $\Omega$ but is not in $\es(\sig)$. Figure~\ref{fig_ext_star_not_contractible} also illustrates the difference between $\es(\sig)$ non-contractible and $\clos(\es(\sig))$ non-contractible.
\end{remark}

% \subsection{Traces on discrete geometric objects}

% Let $l\in\{0{:}3\}$ and let $u \in V_p^l(\Th)$. 
% Let $m\in\{l{:}3\}$ and let $\sigma \in \Delta_h^{m}$ be an oriented geometric object
% associated with the mesh (a vertex, an edge, a face or a tetrahedron). 
% In what follows, we consider the trace $\tr^l_\sigma(u)$
% of $u$ on $\sigma$, defined according to the following table:
% \begin{center}
% definition of $\tr^l_\sigma(u)$, $u\in V_p^l(\Th)$, $\sigma \in \Delta_h^m$, $0\le l\le m\le 3$
% \begin{tabular}{||c | c c c c||} 
%  \hline
%   & $\sig\in\Delta_h^0$ & $\sig\in\Delta_h^1$ & $\sig\in\Delta_h^2$ & $\sig\in\Delta_h^3$ \\ [0.5ex] 
%  \hline\hline
%  $l=0$ & $u(\sig)$ & $u|_\sig$ & $u|_\sig$ & $u|_\sig$\\ 
%  \hline
%  $l=1$ &  & $\bu|_\sig\SCAL \bt_\sig$ & $(\bn_\sig\times \bu|_\sig)\times \bn_\sig$ & $\bu|_\sig$\\
%  \hline
%  $l=2$ &  &  & $\bu|_\sig\SCAL \bn_\sig$ & $\bu|_\sig$\\
%  \hline
%  $l=3$ &  &  &  & $u|_\sig$\\
%  \hline
% \end{tabular}
% \end{center}
% It is well-known that $\tr^l_\sigma(u)$ is well-defined for all $l\in\{0{:}3\}$ and all 
% $u \in V_p^l(\Th)$.

\subsection{Piecewise polynomial spaces, canonical degrees of freedom, and Whitney forms}
\label{sec:Whitney}

Let $p\ge0$ be the polynomial degree. For a tetrahedron $\tau \in \Th$, $\pol_p(\tau)$  is the space of polynomials of degree at most $p$ defined on $\tau$,
$\Ne_p(\tau):=\{\bu(\bx) + \bx\times \bv(\bx) : \bu,\bv\in \pol_p(\tau;\R^3)\}$ is the $p$-th order N\'ed\'elec space~\cite{nedelec1980mixed}, and 
$\Rt_p(\tau):=\{\bu(\bx)+v(\bx)\bx: \bu\in \pol_p(\tau;\R^3), v\in\pol_p(\tau)\}$ is the $p$-th order Raviart--Thomas space~\cite{Ra_Tho_MFE_77}. We consider the following piecewise polynomial spaces:
\begin{subequations} \label{eq_spaces_disc} \begin{align}
    V_p^0(\Th) &:=  \{ u\in H^1(\Omega): u|_\tau \in \pol_{p+1}(\tau), \forall \tau \in \Th\}, \\
    \bV_p^1(\Th) &:=  \{ \bu\in \bH(\curl, \Omega): \bu|_\tau \in \Ne_p(\tau), \forall \tau \in \Th\}, \\
    \bV_p^2(\Th) &:= \{ \bu\in \bH(\dive, \Omega): \bu|_\tau \in \Rt_p(\tau), \forall \tau \in \Th\}, \\
    V_p^3(\Th) &:= \{ u\in L^2(\Omega): u|_\tau \in \pol_{p}(\tau), \forall \tau \in \Th\}.
\end{align} \end{subequations}
We use the generic notation $V_p^l(\Th)$ for the above spaces, and note that $V_p^l(\Th) \subset V^l$, for all $l\in\{0{:}3\}$, with the graph spaces $V^l$ defined in~\eqref{eq:graph_spaces}.

The canonical degrees of freedom are linear forms on $V_p^l(\Th)$, for all $l\in\{0{:}3\}$,
associated with the (oriented) $l$-simplices of the mesh $\Th$. In this work, we only consider explicitly the lowest-order canonical degrees of freedom which are defined as follows:
\begin{subequations}\begin{alignat}{2}
\phi_v(u)&:=u(v),&\quad&\forall u\in V_p^0(\Th),\; \forall v \in \Delta_h^0=\Vh,\\
\phi_e(\bu)&:=\int_e \bu\SCAL\bt_e,&\quad&\forall \bu\in \bV_p^1(\Th),\; \forall e \in \Delta_h^1=\Eh,\\
\phi_f(\bu)&:=\int_f \bu\SCAL\bn_f,&\quad&\forall \bu\in \bV_p^2(\Th),\; \forall f \in \Delta_h^2=\Fh,\\
\phi_\tau(u)&:=\int_\tau u,&\quad&\forall u\in V_p^3(\Th),\; \forall \tau \in \Delta_h^3=\Th.
\end{alignat}\end{subequations}
The above integrals are understood in algebraic form; for instance,
$\int_\tau 1 = |\tau| \Sign(\det(\bt_{e_1},\bt_{e_2},\bt_{e_3}))$ with $\bt_{e_k}$ pointing from $x_{i_0}$ to $x_{i_k}$ for all $k\in\{1{:}3\}$.
It is well-known that $(\phi_v)_{v\in\Vh}$, $(\phi_e)_{e\in\Eh}$, $(\phi_f)_{f\in\Fh}$, $(\phi_\tau)_{\tau\in\Th}$ form a basis for the dual space of the lowest-order piecewise polynomial spaces $V_0^0(\Th)$, $\bV_0^1(\Th)$, $\bV_0^2(\Th)$, $V_0^3(\Th)$, respectively. The higher-order degrees of freedom are also relevant in the construction of the commuting projections, but are not used explicitly in the paper, as this part of the construction follows~\cite{Arn_Guz_loc_stab_L2_com_proj_21}.

The dual bases of the canonical degrees of freedom are composed of the so-called Whitney forms $(W_v)_{v\in\Vh}$, $(\bW_e)_{e\in\Eh}$, $(\bW_f)_{f\in\Fh}$, $(W_\tau)_{\tau\in\Th}$. By construction, the Whitney forms are such that
\begin{equation}
V_0^0(\Th) = \Span_{v\in\Vh} W_v, \quad \bV_0^1(\Th) = \Span_{e\in\Eh} \bW_e, \quad \bV_0^2(\Th) = \Span_{f\in\Fh} \bW_f, \quad V_0^3(\Th) = \Span_{\tau\in\Th} W_\tau,
\end{equation}
and
\begin{equation} \label{eq:dof_Whitney}
\phi_{v'}(W_v)=\delta_{vv'}, \quad
\phi_{e'}(\bW_e)=\delta_{ee'}, \quad
\phi_{f'}(\bW_f)=\delta_{ff'}, \quad
\phi_{\tau'}(W_\tau)=\delta_{\tau\tau'},
\end{equation}
where the $\delta$'s are Kronecker deltas. Moreover, 
for all $l\in\{0{:}3\}$ and all $\sigma\in\Delta_h^l$, 
the support of $W_\sigma$ is $\clos(\st(\sigma))$, and we have (see, e.g., \cite{BBF13, ErnGuermondbook})
\begin{equation}\label{CW}
\|W_\sigma\|_{L^2(\st(\sigma))} \le C_W h_{\sigma}^{\frac32-l} \quad  \forall \sigma \in \Delta_h^l,
\; \forall l\in\{0{:}3\}, 
\end{equation}
where $C_W$ only depends on the shape-regularity parameter $\rho_{\Th}$ of the mesh $\Th$. 

\subsection{Incidence matrices}
\label{sec:incidence}

Incidence matrices are the algebraic realization of the 
differential operators from the de Rham complex~\eqref{complex}, 
but acting on the lowest-order degrees of freedom.
These matrices have emerged in various contexts related to compatible (or mimetic, or
structure-preserving) discretizations. We refer the reader 
to~\cite{Bossa:00,BocHy:05,BonelleErn2014,Gerri:12} and the
references therein for further insight into this topic. 
The incidence matrices of interest here are the compatible discrete gradient,
curl, and divergence matrices such that $\sfG\in \R^{\#\Eh,\#\Vh}$, 
$\sfC\in \R^{\#\Fh,\#\Eh}$, and $\sfD\in \R^{\#\Th,\#\Fh}$, where $\#\bullet$ denotes the
\rev{cardinality} of the finite set $\bullet$.
The incidence matrices satisfy the following complex (or compatibility) properties:
\begin{equation}
\label{eq:compatible}
\sfG \sfU = \sfz, \qquad
\sfC \sfG = \sfz, \qquad
\sfD \sfC = \sfz,
\end{equation}
where $\sfU\in \R^{\#\Vh}$ is the column vector having all of its entries equal to one,
and $\sfz$ is the zero vector or matrix of appropriate size depending on the context.
The above properties are the discrete counterpart
of the following well-known relations satisfied by the differential operators:
$\grad \, c=\bzero$ with $c$ a constant function, $\curl \,\grad=\bzero$, and $\div\curl=0$.

The entries of the incidence matrices are incidence numbers in $\{-1,0,1\}$ associated with 
pairs of oriented geometric objects. To define these
numbers, we introduce the subsets
\begin{subequations} \begin{alignat}{2}
\Ev&:=\{e\in\Eh\,:\, v\in e\},
&&\quad \forall v\in\Vh, \\
\Fe&:=\{f\in\Fh\,:\, e\subset f\}, &&\quad \forall e\in\Eh, \\
\Tf&:=\{\tau\in\Th\,:\, f\subset \tau\}, &&\quad 
\forall f\in\Fh.
\end{alignat} \end{subequations}
If $e:=[x_{i_0},x_{i_1}] \in \Ev$, $v$ is obtained from $e$ by omitting one of the two vertices
of $e$, say $x_{i_j}$ with $j\in\{0,1\}$, and we set $\iota_{ev}:=(-1)^j$. 
If $f:=[x_{i_0},x_{i_1},x_{i_2}] \in \Fe$, $e$ is obtained from $f$ by omitting one of the three vertices of $f$, say $x_{i_j}$ with $j\in\{0,1,2\}$, and we set $\iota_{fe}:=(-1)^j$. 
If $\tau=[x_{i_0},x_{i_1},x_{i_2},x_{i_3}] \in \Tf$, $f$ is obtained from $\tau$ by omitting one of the four vertices of $\tau$, say $x_{i_j}$ with $j\in\{0,1,2,3\}$, and we set $\iota_{\tau f}:=(-1)^j$.
Finally, we also set $\iota_{ev}:=0$ for all $v\in \Vh$ and all $e\not\in\Ev$,
$\iota_{fe}:=0$ for all $e\in \Eh$ and all $f\not\in\Fe$, and
$\iota_{\tau f}:=0$ for all $f\in \Fh$ and all $\tau\not\in\Tf$.

The incidence matrices have entries such that
\begin{equation}
\sfG_{ev}=\iota_{ev},\; \forall (e,v)\in \Eh{\times}\Vh, \quad
\sfC_{fe}=\iota_{fe},\; \forall (f,e)\in \Fh{\times}\Eh, \quad
\sfD_{\tau f}=\iota_{\tau f},\; \forall (\tau,f)\in \Th{\times}\Fh.
\end{equation} 
The complex properties~\eqref{eq:compatible} take the following form:
\begin{subequations} \label{eq:compatible_iota} \begin{alignat}{2}
\sum_{v\in\Ve} \iota_{ev} &=0, \quad &&\forall e\in\Eh, \label{eq:GRAD} \\
\sum_{e\in\Ev\cap \Ef} \iota_{fe}\iota_{ev} &=0, \quad && \forall f\in \Fh, \; \forall v\in \Vh, \label{eq:CURL_GRAD} \\
\sum_{f\in\Fe \cap \Ft} \iota_{\tau f}\iota_{fe} &=0, \quad && \forall \tau\in \Th, \; \forall e\in \Eh. \label{eq:DIV_CURL}
\end{alignat} \end{subequations}

The following identities are a straightforward consequence of the Stokes theorem~\cite{whitney2012geometric}:
\begin{subequations} \begin{alignat}{3}
\phi_e(\grad \, m) = \sum_{v\in \Ve} \iota_{ev} \phi_v(m), \quad &&\forall e\in\Eh, \quad &&\forall m\in V^0_p(\Th), \label{eq:dofs_iota_0}\\
\phi_f(\curl \, \bbm) = \sum_{e\in \Ef} \iota_{fe} \phi_e(\bbm), \quad &&\forall f\in\Fh, \quad &&\forall \bbm\in \bV^1_p(\Th), \label{eq:dofs_iota_1}\\
\phi_\tau(\dive \, \bbm) = \sum_{f\in \Ft} \iota_{\tau f} \phi_f(\bbm), \quad &&\forall \tau\in\Th, \quad &&\forall \bbm\in \bV^2_p(\Th). \label{eq:dofs_iota_2}
\end{alignat} \end{subequations}
Notice that the above identities are valid more generally for smooth functions and fields. Applying the above identities to the Whitney forms, one readily infers that
\begin{equation} \label{eq:W_iota}
\grad \, W_v = \sum_{e\in \Ev} \iota_{ev} \bW_e,
\qquad
\curl \, \bW_e = \sum_{f\in \Fe} \iota_{fe} \bW_f,
\qquad
\div \bW_f = \sum_{\tau\in \Tf} \iota_{\tau f} W_\tau.
\end{equation} 

\begin{remark}[Link to FEEC formalism]
Incidence matrices are not explicitly present in the FEEC formalism, but are replaced 
by the notion of boundary operators acting on $l$-chains, which are formal (algebraic) sums 
of $l$-simplices in $\Th$ for all $l\in\{0{:}3\}$. For an oriented $l$-simplex 
$\sigma := [x_{i_0},\ldots,x_{i_l}]$ and a permutation $\varphi:\{0{:}l\}\to \{0{:}l\}$,
one can define the geometric object $\sigma^\varphi:=[x_{i_{\varphi(0)}},\ldots,x_{i_{\varphi(l)}}]$.
Notice that $\sigma$ and $\sigma^\varphi$ coincide as sets of points in $\R^3$, but may have
a different orientation. We write $\sigma^\varphi=\sigma$ if $\varphi$ is an even permutation (both objects then have the same orientation), and $\sigma^\varphi=-\sigma$ if $\varphi$ is an odd permutation (the two objects then have an opposite orientation). Then, letting $\calC^l$ be the space composed of formal sums of $l$-simplices, the boundary operator $\partial_l : \calC^l\rightarrow \calC^{l-1}$, for all $l\in\{1{:}3\}$, is the linear operator such that, for all $\sigma := [x_{i_0},\ldots,x_{i_l}]$,
\[
\partial_l \sigma := \sum_{j\in\{0{:}l\}} (-1)^j [x_{i_0},\ldots,\hat x_{i_j},\ldots,x_{i_l}],
\]
where $\hat x_{i_j}$ means omission of $x_{i_j}$. The links between the boundary operator and the incidence matrices are as follows:
\begin{alignat*}{2}
\partial_1 e &= \sum_{v\in \Ve} \iota_{ev} v, \quad &&\forall e\in \Eh, \\
\partial_2 f &= \sum_{e\in \Ef} \iota_{fe} e, \quad &&\forall f\in \Fh, \\
\partial_3 \tau &= \sum_{f\in \Ft} \iota_{\tau f} f, \quad &&\forall \tau\in \Th.
\end{alignat*}
Moreover, the counterpart of the complex properties~\eqref{eq:GRAD}--\eqref{eq:DIV_CURL} are
\[
\partial_0 \partial_1 e = 0 \; \forall e\in \Eh,\qquad
\partial_1 \partial_2 f = 0 \; \forall f\in \Fh,\qquad
\partial_2 \partial_3 \tau = 0 \; \forall \tau \in \Th,
\]
with the convention that $\partial_0v = 1$ for all $v\in\Vh$.
\end{remark}

\section{Local weight functions}
\label{sec:weight}

The goal of this section is to state our main result concerning the local weight functions.
\rev{In contrast to~\cite{Arn_Guz_loc_stab_L2_com_proj_21}, we build local weight functions that are discrete (piecewise polynomials). A simple and convenient way to do so is to employ the so-called Alfeld split of the mesh $\Th$ that we recall below. Piecewise polynomial weight functions then allow us to invoke discrete Poincar\'e inequalities (as opposed to continuous right-inverses in $H^1$) to establish the $L^2$-stability of the weight functions used to define the projections and thus to bypass~\cite[Assumption~2.4]{Arn_Guz_loc_stab_L2_com_proj_21}.
Using the Alfeld split allows us to work with piecewise \emph{affine} bubble functions and, thus, consider the \emph{same} polynomial degree for the discrete Poincar\'e inequalities on the Alfeld split. Alternatively, one could consider the given mesh only, together with the standard higher-order (4th-degree for $d=3$) bubble functions. This would increase the polynomial degree needed for the local weight functions of the starting space $V^3_p$ and also on each consecutive de Rham step since the local weight functions are constructed progressively (backwards).}
To streamline the presentation, the proofs of the results stated in this section are given in the next section and in the appendix.

\subsection{Alfeld split}

\begin{figure}[H]
  \centering
  \definecolor{grayish}{rgb}{0.36, 0.54, 0.66}
  \begin{subfigure}[b]{0.4\textwidth}
     \centering
         \begin{tikzpicture}[scale=0.66]
         % Define vertices
         \coordinate (x0) at (4,0);
         \coordinate (x1) at (0,0);
         \coordinate (x2) at (2,3.5);
         \coordinate (x3) at (1.5,1.0);

         % Draw tetrahedron faces
         \fill[grayish!50,opacity=0.7] (x0) -- (x1) -- (x3) -- cycle;  % Base triangle x0, x1, x3
         \fill[grayish!30,opacity=0.8] (x0) -- (x2) -- (x3) -- cycle;  % Side face x0, x2, x3
         \fill[grayish!40,opacity=0.9] (x1) -- (x2) -- (x3) -- cycle;  % Side face x1, x2, x3
         \fill[grayish!60,opacity=0.5] (x0) -- (x1) -- (x2) -- cycle;  % Top face x0, x1, x2
         
         % Draw outer edges
         \draw[thick] (x0) -- (x1) -- (x2) -- cycle;  % Top triangle x0, x1, x2
         \draw[thick] (x0) -- (x3);  % Edge x0-x3
         \draw[thick] (x1) -- (x3);  % Edge x1-x3
         \draw[thick] (x2) -- (x3);  % Edge x2-x3
         
         % Label vertices
         \node[below right] at (x0) {$x_0$};
         \node[below left] at (x1) {$x_1$};
         \node[above] at (x2) {$x_2$};
         \node[above right] at (x3) {$x_3$};
         
         \end{tikzpicture}
     \caption{Single tetrahedron $\tau$}
  \end{subfigure}
  \hspace{2mm}
  \begin{subfigure}[b]{0.4\textwidth}
     \centering
         \begin{tikzpicture}[scale=0.66]
         % Define vertices
         \coordinate (x0) at (4,0);
         \coordinate (x1) at (0,0);
         \coordinate (x2) at (2,3.5);
         \coordinate (x3) at (1.5,1.0);
         
         % Define the barycenter (z) as the average of the four vertices
         \coordinate (z) at ($(x0)!0.25!(x1)!0.25!(x2)!0.25!(x3)$);
         
         % Draw tetrahedron faces
         \fill[grayish!50,opacity=0.7] (x0) -- (x1) -- (x3) -- cycle;  % Base triangle x0, x1, x3
         \fill[grayish!30,opacity=0.8] (x0) -- (x2) -- (x3) -- cycle;  % Side face x0, x2, x3
         \fill[grayish!40,opacity=0.9] (x1) -- (x2) -- (x3) -- cycle;  % Side face x1, x2, x3
         \fill[grayish!60,opacity=0.5] (x0) -- (x1) -- (x2) -- cycle;  % Top face x0, x1, x2
         
         % Draw outer edges
         \draw[thick] (x0) -- (x1) -- (x2) -- cycle;  % Top triangle x0, x1, x2
         \draw[thick] (x0) -- (x3);  % Edge x0-x3
         \draw[thick] (x1) -- (x3);  % Edge x1-x3
         \draw[thick] (x2) -- (x3);  % Edge x2-x3
         
         % Draw internal dashed lines connecting to the barycenter z
         \draw[dashed, thick, black] (x0) -- (z);
         \draw[dashed, thick, black] (x1) -- (z);
         \draw[dashed, thick, black] (x2) -- (z);
         \draw[dashed, thick, black] (x3) -- (z);

         % Label vertices
         \node[below right] at (x0) {$x_0$};
         \node[below left] at (x1) {$x_1$};
         \node[above] at (x2) {$x_2$};
         \node[above right] at (x3) {$x_3$};
         \node[above right] at (z) {$\overline{x}$};  % Label the barycenter
         
          % Draw nodes for vertices
        \foreach \point in {x0,x1,x2,x3,z}
          \fill[black] (\point) circle (2pt);

         \end{tikzpicture}
     \caption{Alfeld split of $\tau$}
  \end{subfigure}
     \caption{Example of an Alfeld split}
     \label{fig:Alfeld}
\end{figure}
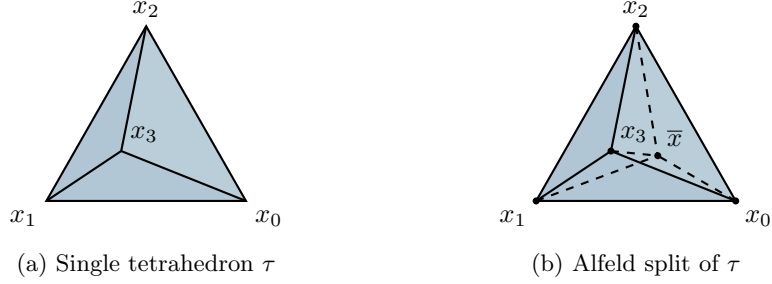

To build the local weight functions, it will be helpful to consider the Alfeld split of $\Th$~\cite{Alfel84}. For every tetrahedron $\tau\in\Th$, we add the barycenter of $\tau$ and connect it to the vertices of $\tau$. This produces four new tetrahedra as illustrated in Figure~\ref{fig:Alfeld}.  We call the resulting mesh $\ThA$. 

% \begin{figure}[H]
%      \centering
%      \begin{subfigure}[b]{0.4\textwidth}
%         \centering
%         \includegraphics[width = 0.4\textwidth]{figures/triangulation.png}
%         \caption{Triangulation $\Th$}
%      \end{subfigure}
%      \hspace{2mm}
%      \begin{subfigure}[b]{0.4\textwidth}
%         \centering
%         \includegraphics[width = 0.4\textwidth]{figures/alfeld.png}
%         \caption{Alfeld split $\ThA$}
%      \end{subfigure}
%         \label{fig:subtri}
%         \caption{Example of Alfeld split in two dimensions.}
% \end{figure}

We let $\mu$ be the globally continuous function  that is piecewise affine on the mesh $\ThA$ such that it vanishes on $\partial \tau$ for every $\tau \in \Th$ and takes the value one at its barycenter. 
For all $\sigma \in \Delta_h$, we let 
\begin{equation} \label{eq:def_mu}
\mu_{\sigma} := \chi_{\es(\sigma)} \mu,
\end{equation}
where $\chi_{\es(\sigma)}$ is the characteristic function of $\es(\sigma)$.

\subsection{Local polynomial spaces and discrete Poincar\'e inequalities}

Recall that $p\ge0$ denotes the polynomial degree.
The piecewise polynomial spaces $V_p^l(\Th)$ defined in~\eqref{eq_spaces_disc}
lead to the following sub-complex of~\eqref{complex}:
\begin{alignat}{4}\label{discretecomplex}
&\mathbb{R}
\stackrel{\subset}{\xrightarrow{\hspace*{0.5cm}}}\
 V_p^0(\Th)
&&\stackrel{\grad}{\xrightarrow{\hspace*{0.5cm}}}\
 \bV_p^1(\Th)
&&\stackrel{\curl}{\xrightarrow{\hspace*{0.5cm}}}\
 \bV_p^2(\Th)
&&\stackrel{\dive}{\xrightarrow{\hspace*{0.5cm}}}\
V_p^3(\Th)
\stackrel{}{\xrightarrow{\hspace*{0.5cm}}}\
0.
\end{alignat}

In the construction of the local weight functions, we need to consider local
versions of these spaces. 
For all $l\in\{0{:}3\}$ and all $\sigma \in \Delta_h$ 
(notice that $\sigma$ is not necessarily in $\Delta_h^l$), 
we let $\calTes$ be the collection of tetrahedra in $\Th$ composing $\es(\sigma)$
and let $V_p^l(\calTes)$ be the restriction of the space $V_p^l(\Th)$ to $\es(\sigma)$, i.e.,
\begin{equation} \label{eq:local_spaces}
V_p^l(\calTes):=\{ u|_{\es(\sigma)}\,:\, u\in V_p^l(\Th)\}.
\end{equation}  
We also need local restrictions of the polynomial spaces $V_p^l(\ThA)$ which are
built as the piecewise polynomial spaces $V_p^l(\Th)$ in~\eqref{eq_spaces_disc}, but on the Alfeld split $\ThA$. 
For these local spaces, we enforce homogeneous boundary conditions on the boundary $\partial\es(\sigma)$ of the extended star $\es(\sigma)$ or zero mean-value on $\es(\sigma)$ as follows:
Letting $\calTesA$ be the collection of the split tetrahedra in $\ThA$ composing $\es(\sigma)$, we define
\begin{subequations} \label{eq:local_spaces_m} \begin{align}
\mVp{0}(\calTesA) &:= \{ u|_{\es(\sigma)}\,:\, u \in V_p^0(\ThA),\, u|_{\partial \es(\sigma)}=0\}, \\
\mbVp{1}(\calTesA) &:= \{ \bu|_{\es(\sigma)}\,:\, \bu \in \bV_p^1(\ThA), \, \bu|_{\partial \es(\sigma)}{\times}\bn_{\es(\sigma)}=\bzero \}, \\
\mbVp{2}(\calTesA) &:= \{ \bu|_{\es(\sigma)}\,:\, \bu \in \bV_p^2(\ThA),\, \bu|_{\partial \es(\sigma)}{\cdot}\bn_{\es(\sigma)}=0 \}, \\
\mVp{3}(\calTesA) &:= \{ u|_{\es(\sigma)}\,:\, u \in V_p^3(\ThA), \, \bl u,1\br_{\es(\sigma)}=0\}, 
\end{align} \end{subequations}
where $\bn_{\es(\sigma)}$ denotes the unit outward normal to $\es(\sigma)$. 
We generically denote the above spaces as $\mVp{l}(\calTesA)$ for all $l\in\{0{:}3\}$.
We also define $V_p^3(\calTesA) := \{ u|_{\es(\sigma)}\,:\, u \in V_p^3(\ThA)\}$.

We recall the following result, which is noted, for instance, in \cite[Section 2.2]{Falk_Winth_loc_coch_14}.

\begin{proposition}[Exact sequences]\label{exactcomplex}
For all $\sigma \in \Delta_h$, under the assumption that $\clos(\es(\sigma))$ is contractible, the following discrete sequences are exact: 
\begin{subequations} 
\begin{equation}\label{localcomplex}
\mathbb{R}
\stackrel{\subset}{\xrightarrow{\hspace*{0.5cm}}}\
 V_p^0(\calTes)
\stackrel{\grad}{\xrightarrow{\hspace*{0.5cm}}}\
  \bV_p^1(\calTes)
\stackrel{\curl}{\xrightarrow{\hspace*{0.5cm}}}\
  \bV_p^2(\calTes)
\stackrel{\dive}{\xrightarrow{\hspace*{0.5cm}}}\
 V_p^3(\calTes)
\stackrel{}{\xrightarrow{\hspace*{0.5cm}}}\
0,
\end{equation}
\begin{equation}\label{localcomplexm}
0 
\stackrel{\subset}{\xrightarrow{\hspace*{0.5cm}}}\
 \mVp{0}(\calTesA)
\stackrel{\grad}{\xrightarrow{\hspace*{0.5cm}}}\
  \mbVp{1} (\calTesA)
\stackrel{\curl}{\xrightarrow{\hspace*{0.5cm}}}\
  \mbVp{2}(\calTesA)
\stackrel{\dive}{\xrightarrow{\hspace*{0.5cm}}}\
 \mVp{3}(\calTesA)
\stackrel{}{\xrightarrow{\hspace*{0.5cm}}}\
0.
\end{equation}
\end{subequations}
\end{proposition}

We define the kernels of the local spaces as 
\begin{subequations} \label{eq:def_ZZ} \begin{align}
\Zz   V_p^0(\calTes)&:=  \{ u \in  V_p^0(\calTes): \grad \, u=\bzero \},  \\
\Zz   \bV_p^1(\calTes)&:= \{ \bu \in  \bV_p^1(\calTes): \curl \, \bu=\bzero \},  \\
\Zz   \bV_p^2(\calTes)&:= \{ \bu \in  \bV_p^2(\calTes): \dive \, \bu=0 \}, 
\end{align} \end{subequations}
and the orthogonal complements as
\begin{subequations} \label{eq:def_ZZ_perp} \begin{align}
\Zz^{\perp}   V_p^0(\calTes)&:= \{ u \in  V_p^0(\calTes): \bl u,v \br_{\es(\sigma)} =0, \forall v \in  \Zz   V_p^0(\calTes)\}, \\ 
\Zz^{\perp}   \bV_p^1(\calTes)&:= \{ \bu \in  \bV_p^1(\calTes): \bl \bu,\bv \br_{\es(\sigma)} =0, \forall \bv \in  \Zz   \bV_p^1(\calTes) \}, \\ 
\Zz^{\perp}   \bV_p^2(\calTes)&:= \{ \bu \in  \bV_p^2(\calTes): \bl \bu,\bv \br_{\es(\sigma)} =0, \forall \bv \in  \Zz   \bV_p^2(\calTes) \}. 
\end{align} \end{subequations}
\rev{Since these spaces are composed of piecewise polynomials, we will be able to establish the following discrete Poincar\'e inequalities on extended stars of arbitrary shape. This represents a conceptual difference with respect to~\cite{Arn_Guz_loc_stab_L2_com_proj_21}, where continuous right-inverses in $H^1$ were needed, but merely supposed in Assumption~2.4 on non-star-shaped extended stars $\es(\sigma)$.}

\begin{proposition}[Discrete Poincar\'e inequalities on extended stars] \label{prop:discP_star}
There exists a constant $\CP$, only depending on the mesh shape-regularity parameter $\rho_{\Th}$ and the polynomial degree $p$, such that, for all $\sigma \in \Delta_h$,
\begin{subequations}
\begin{alignat}{2}
\|u\|_{L^2(\es(\sigma))} &\le \CP h_{\sigma} \| \grad \, u \|_{\bL^2(\es(\sigma))},  \qquad && \forall u \in \Zz^{\perp} V_p^0(\calTes),   \label{onto0}\\
\|\bu\|_{\bL^2(\es(\sigma))} &\le  \CP h_{\sigma} \|\curl \, \bu\|_{\bL^2(\es(\sigma))}, \qquad && \forall  \bu \in \Zz^{\perp} \bV_p^1(\calTes), \label{onto1} \\
\|\bu\|_{\bL^2(\es(\sigma))} &\le  \CP h_{\sigma} \|\dive \, \bu\|_{L^2(\es(\sigma))},   \qquad && \forall  \bu \in \Zz^{\perp} \bV_p^2(\calTes).  \label{onto2}
\end{alignat}
\end{subequations}
\end{proposition}

We define the spaces  $\Zz \mVp{l}(\calTesA)$ and $\Zz^{\perp} \mVp{l}(\calTesA)$ as in~\eqref{eq:def_ZZ}--\eqref{eq:def_ZZ_perp}. Similarly to Proposition~\ref{prop:discP_star}, we will establish the following discrete Poincar\'e inequalities with boundary conditions on Alfeld splits of extended stars.

\begin{proposition}[Discrete Poincar\'e inequalities on Alfeld splits of extended stars] \label{prop:discP_star_Alfeld}
There exists a constant $\CPA$, only depending on the mesh shape-regularity parameter $\rho_{\Th}$ and the polynomial degree $p$,  such that, for all $\sigma \in \Delta_h$,
\begin{subequations}
\begin{alignat}{2}
 \|u\|_{L^2(\es(\sigma))} &\le \CPA h_{\sigma} \| \grad \, u \|_{\bL^2(\es(\sigma))},  \qquad && \forall u \in \Zz^{\perp} \mVp{0}(\calTesA),   \label{onto0m}\\
 \|\bu\|_{\bL^2(\es(\sigma))} &\le  \CPA h_{\sigma} \|\curl \, \bu\|_{\bL^2(\es(\sigma))}, \qquad && \forall  \bu \in \Zz^{\perp} \mbVp{1}(\calTesA), \label{onto1m} \\
 \|\bu\|_{\bL^2(\es(\sigma))} &\le  \CPA h_{\sigma} \|\dive \, \bu\|_{L^2(\es(\sigma))},   \qquad && \forall  \bu \in \Zz^{\perp} \mbVp{2}(\calTesA).  \label{onto2m}
\end{alignat}
\end{subequations}
\end{proposition}
The proof of Propositions~\ref{prop:discP_star} and~\ref{prop:discP_star_Alfeld} is postponed to Appendix~\ref{sec:proof_Poinc}.

\subsection{Main result}

We are now ready to state our main result. The proof is postponed to Section~\ref{sec:proof}.

%There exists linear operators $\ZZ_p^0: \C_0 \mapsto \mV_p^{3,A},  \ZZ_p^1:\C_1   \mapsto \mV_p^{2,A} , \ZZ_p^2:\C_2 \mapsto  \mV_p^{1,A} , \ZZ_p^3: \C_3\mapsto   \mV_p^{0,A}$ 

\begin{theorem}[Local weight functions] \label{th:weights}
Let $p\ge0$ be the polynomial degree. There exist local weight functions
%a family of local weight functions $\{\ZZ_p^l(\sigma)\}_{\sigma \in \Delta_h^l}$ with $\ZZ_p^l(\sigma) \in \mVp{3-l}(\es(\sigma))$, that is,
$\ZZ_p^0(v) \in V_p^3(\calTevA)$ for all vertices $v \in \Vh$ (piecewise polynomial of order $p$ on the Alfeld split of the (extended) vertex star $\es(v)$, no mean-value condition),
$\bZZ_p^1(e) \in \mbVp{2}(\calTeeA)$ for all edges $e\in \Eh$ (piecewise Raviart--Thomas polynomial of order $p$ on the Alfeld split of the extended edge star $\es(e)$ with zero normal component on $\partial \es(e)$),
$\bZZ_p^2(f) \in \mbVp{1}(\calTefA)$ for all faces $f\in \Fh$ (piecewise N\'ed\'elec polynomial of order $p$ on the Alfeld split of the extended face star $\es(f)$ with zero tangential component on $\partial \es(f)$), and
$\ZZ_p^3(\tau) \in \mVp{0}(\calTetA)$ for all elements $\tau \in\Th$ (continuous piecewise  polynomial of degree $(p+1)$ on the Alfeld split of the extended element star $\es(\tau)$ with zero value on $\partial \es(\tau)$),
satisfying the following properties:\\
\textup{(i)} Support and $L^2$-norm: For all $l\in\{0{:}3\}$ and all $\sig \in \Delta_h^l$,
\begin{equation}\label{ZZc}
\supp\ZZ_p^l(\sig) \subseteq \clos(\es(\sig)),  \qquad
\|\ZZ_p^l(\sig)\|_{L^2(\es(\sig))}  \le C_Z h_{\sig} ^{-\frac{3}{2}+l},
\end{equation}
where $C_Z$ only depends on the mesh shape-regularity parameter $\rho_{\Th}$ and the polynomial degree $p$.
\\
\textup{(ii)} Relation to canonical degrees of freedom:
\begin{subequations}\label{ZZb}
\begin{alignat}{3}
 \bl \ZZ_p^0(v), u \br_{\es(v)}&=  \phi_v(u), \qquad && \forall u \in V_p^0(\Th), &&\quad \forall v \in  \Delta_h^0=\Vh, \label{ZZb0}\\
 \bl \bZZ_p^1(e), \bu \br_{\es(e)}&=  \phi_e(\bu), \qquad && \forall \bu \in \bV_p^1(\Th), &&\quad \forall e \in  \Delta_h^1=\Eh, \label{ZZb1} \\
  \bl \bZZ_p^2(f), \bu \br_{\es(f)}&=  \phi_f(\bu), \qquad && \forall \bu \in \bV_p^2(\Th), &&\quad \forall f \in  \Delta_h^2=\Fh,  \label{ZZb2}\\
 \bl \ZZ_p^3(\tau), u \br_{\es(\tau)}&=  \phi_\tau(u),  \qquad && \forall u \in V_p^3(\Th), &&\quad \forall \tau \in  \Delta_h^3=\Th. \label{ZZb3}
\end{alignat}
\end{subequations}
\textup{(iii)} Relation to differential operators:
\begin{subequations}\label{ZZa}
\begin{alignat}{2}
-\dive \, \bZZ_p^1(e)&= \sum_{v\in \Ve} \iota_{ev} \ZZ_p^0(v), \qquad && \forall e \in \Delta_h^1=\Eh,  \label{ZZa0}\\
\curl \, \bZZ_p^2(f)&= \sum_{e\in \Ef} \iota_{fe} \bZZ_p^1(e), \qquad && \forall f \in \Delta_h^2=\Fh, \label{ZZa1} \\
-\grad \, \ZZ_p^3(\tau)&= \sum_{f\in \Ft} \iota_{\tau f} \bZZ_p^2(f), \qquad && \forall \tau \in \Delta_h^3=\Th, \label{ZZa2}
\end{alignat} \end{subequations}
with the incidence numbers defined in Section~\textup{\ref{sec:incidence}}.
\end{theorem}

These local weight functions are the key ingredient to construct the local $L^2$-bounded commuting projections. This is done in Section~\ref{sec:projections}. 
% In particular, the lowest-order ($p=0$) projection only requires the lowest-order weight functions.

\section{Proof of Theorem~\ref{th:weights} on the local weight functions}
\label{sec:proof}

In this section, we prove our main result on the weight functions (Theorem~\ref{th:weights}).
The proof relies on a constructive argument; indeed, we construct
explicitly the weight functions satisfying the expected properties. 
For positive numbers $a,b$, we abbreviate as $a\lesssim b$ the
inequality $a\le Cb$, where the value of the positive constant $C$ can change at each
occurrence, but it can only depend on the shape-regularity parameter $\rho_{\Th}$ of the mesh $\Th$ and the polynomial degree $p$. 
We write $a\approx b$ when both inequalities $a\lesssim b$ and $b\lesssim a$ hold.

\subsection{Construction and properties of $\ZZ_p^0(v)$}

Let $v \in \mathcal{V}_h$ and let 
\begin{equation} \label{eq_eta_0} 
    \eta_0(v):=\chi_{\es(v)}/|\es(v)|,
\end{equation} 
where we recall that $\chi_{\es(v)}$ is the characteristic function of $\es(v)$. 
We define $\psi_0(v) \in \Zz^{\perp}  V_p^0(\calTev)$ (a continuous piecewise polynomial of degree $(p+1)$ on the original mesh of the (extended) vertex star $\es(v)$ with zero mean-value on $\es(v)$) by
\begin{equation}\label{605}
\bl \mu_{v} \, \grad \, \psi_0(v), \grad \, u\br_{\es(v)} = \phi_v(u) -\bl \eta_0(v), u \br_{\es(v)}, \quad  \forall u \in   \Zz^{\perp} V_p^0(\calTev),
\end{equation}
with $\mu_v$ defined in~\eqref{eq:def_mu}.
This problem is well-posed owing to the Poincar\'e inequality~\eqref{onto0}. 
Notice that the right-hand side of~\eqref{605} is given by the point value of the test function $u$ at the vertex $v$ minus the mean-value of $u$ on the (extended) vertex star $\es(v)$. Thus,
both the left-hand side and the right-hand side of~\eqref{605} vanish when $u\in 
\Zz V_p^0(\calTev)$ (i.e., on constant functions), and we infer that 
\begin{equation}\label{805}
\bl \mu_{v} \, \grad \, \psi_0(v), \grad \, u\br_{\es(v)}=\phi_v(u) -\bl \eta_0(v), u \br_{\es(v)}, \quad  \forall u \in V_p^0(\calTev).
\end{equation}

We now define
\begin{equation}\label{defZ0}
\boxed{\ZZ_p^0(v):=\eta_0(v)- \dive \big(\mu_v \, \grad \, \psi_0(v) \big) \quad \text{in $\es(v)$},}
\end{equation}
and extend $\ZZ_p^0(v)$ by zero outside $\es(v)$. 
Notice that each component of $\grad \, \psi_0(v)$ is a piecewise polynomial of degree $p$ on the original mesh $\Th$ whose normal component is possibly discontinuous on the faces of $\Th$. However, $\mu_v$ is zero on these faces, so $\mu_v \grad \, \psi_0(v)$ is a continuous piecewise polynomial of degree $(p+1)$ on the Alfeld split $\ThA$ and is zero on $\partial(\es(v))$. Thus, we can take the divergence of $\mu_v \grad \, \psi_0(v)$ and obtain a piecewise  polynomial of degree $p$ on the Alfeld split $\calTevA$. Hence, we indeed have $\ZZ_p^0(v)\in V_p^3(\calTevA)$.
It remains to prove~\eqref{ZZc} for $l=0$ and~\eqref{ZZb0}. 

(1) Proof of~\eqref{ZZc} for $l=0$.
By construction, 
$\ZZ_p^0(v)$  has support in $\clos(\es(v))$.
We first observe that 
\begin{equation}\label{eta0bound}
\| \eta_0(v)\|_{L^2(\es(v))} = |\es(v)|^{-\frac12} \lesssim h_v^{-\frac32}. 
\end{equation}
Moreover, using~\eqref{605} and the Cauchy--Schwarz inequality, we obtain 
\begin{alignat*}{2}
\| \grad \, \psi_0(v)\|_{\bL^2(\es(v))}^2 &\lesssim \bl \mu_{v} \, \grad \, \psi_0(v), \grad \, \psi_0(v)\br_{\es(v)} \\
&=\phi_v(\psi_0(v)) -\bl \eta_0(v), \psi_0(v) \br_{\es(v)} \\
&\le \|\psi_0(v)\|_{L^\infty(\es(v))} + \| \eta_0(v)\|_{L^2(\es(v))} \|\psi_0(v)\|_{L^2(\es(v))}.
\end{alignat*}
Hence, invoking the inverse inequality 
\begin{equation} \label{eq_inv_ineq}
    \|\psi_0(v)\|_{L^\infty(\es(v))} \lesssim h_v^{-\frac32} \|\psi_0(v)\|_{L^2(\es(v))}
\end{equation}
together with the Poincar\'e inequality~\eqref{onto0}
to bound $\|\psi_0(v)\|_{L^2(\es(v))}$, we infer that 
\begin{align*}
\| \grad \, \psi_0(v)\|_{\bL^2(\es(v))}^2 &\lesssim (h_v^{-\frac32}+\| \eta_0(v)\|_{L^2(\es(v))})
h_v \| \grad \, \psi_0(v)\|_{\bL^2(\es(v))} \\
&\lesssim h_v^{-\frac32+1}\| \grad \, \psi_0(v)\|_{\bL^2(\es(v))}.
\end{align*}
This proves that
\[
\| \grad \, \psi_0(v)\|_{\bL^2(\es(v))} \lesssim h_v^{-\frac32+1}.
\]
Finally, invoking an inverse estimate, we obtain
\begin{equation*}
  \| \dive \big(\mu_v \, \grad \, \psi_0(v) \big)\| _{L^2(\es(v))}\lesssim h_v^{-1} \| \mu_v \, \grad \, \psi_0(v)\|_{\bL^2(\es(v))} \leq h_v^{-1} \| \grad \, \psi_0(v)\|_{\bL^2(\es(v))} \lesssim h_v^{-\frac32}.
\end{equation*}
Combining the above bounds proves~\eqref{ZZc} for $l=0$.
  
(2) Proof of~\eqref{ZZb0}. Using integration by parts shows that, for all $u \in V_p^0(\calTev)$,
\begin{align*}
\bl \ZZ_p^0(v), u \br_{\es(v)}&= \bl \eta_0(v)- \dive\big( \mu_v \, \grad \, \psi_0(v)\big), u \br_{\es(v)} \\
&= \bl \eta_0 (v), u \br_{\es(v)}+ \bl  \mu_v \, \grad \, \psi_0(v), \grad \,  u \br_{\es(v)} = \phi_v(u),
\end{align*}
owing to~\eqref{805}. Since $\ZZ_p^0(v)$ is zero outside $\es(v)$, the above identity holds
for all $u\in V_p^0(\Th)$. This proves~\eqref{ZZb0}. Notice in passing that the above identity implies that
\begin{equation} \label{eq:int_Z0}
\int_{\es(v)} \ZZ_p^0(v) = 1.
\end{equation}

\subsection{Construction and properties of $\bZZ_p^1(e)$}\label{sec:ZZ1}
Let $e \in \mathcal{E}_h$. Since $\es(v)\subset\es(e)$ for all $v\in\Ve$ and $\ZZ_p^0(v)$ is supported on $\clos(\es(v))$, we infer from~\eqref{eq:int_Z0} that
\begin{equation*}
\int_{\es(e)} \sum_{v\in\Ve} \iota_{ev} \ZZ_p^0(v)
= \sum_{v\in\Ve} \iota_{ev} \int_{\es(v)} \ZZ_p^0(v)= \sum_{v\in\Ve} \iota_{ev} =0,
\end{equation*}
where the last equality follows from~\eqref{eq:GRAD}.
Hence, by the exactness of the discrete sequence~\eqref{localcomplexm} from Proposition~\ref{exactcomplex}, there exists a piecewise Raviart--Thomas polynomial of order $p$ with zero normal component, $\tilde \beeta_1(e) \in \mbVp{2}(\calTeeA)$, such that 
\begin{equation} \label{eq_tilde_beeta}
 -\dive \, \tilde \beeta_1(e)= \sum_{v\in\Ve} \iota_{ev} \ZZ_p^0(v) \quad \text{ on $\es(e)$}.
\end{equation}
Defining $\rev{\tilde \beeta_1^\perp}(e)$ to be the orthogonal projection of $\tilde \beeta_1(e)$ onto $\Zz \mbVp{2}(\calTeeA)$, i.e., $\rev{\tilde \beeta_1^\perp}(e) \in \Zz \mbVp{2}(\calTeeA)$ is the divergence-free Raviart--Thomas piecewise polynomial such that 
\[
\bl \rev{\tilde \beeta_1^\perp}(e), \bu \br_{\es(e)} = \bl \tilde \beeta_1(e), \bu \br_{\es(e)}, \quad  \forall \bu \in \Zz \mbVp{2}(\calTeeA),
\]
we obtain $\beeta_1(e) := \tilde \beeta_1(e) - \rev{\tilde \beeta_1^\perp}(e) \in \Zz^{\perp} \mbVp{2}(\calTeeA)$ such that 
\begin{equation}\label{eta1def}
 -\dive \, \beeta_1(e)= \sum_{v\in\Ve} \iota_{ev} \ZZ_p^0(v) \quad \text{ on $\es(e)$}.
\end{equation} 
% Actually, $\beeta_1(e)$ is the unique solution to the following constrained quadratic minimization problem:
% \begin{equation} \label{eq_RT_min}
%     \beeta_1(e) = \Argmin_{\substack{\bv \in \mbVp{2}(\calTeeA) \\ -\dive \, \beeta_1(e)= \sum_{v\in\Ve} \iota_{ev} \ZZ_p^0(v)}} \|\bv\|_{\bL^2(\es(e))}^2.
% \end{equation}
\rev{From} the Poincar\'e inequality~\eqref{onto2m}, we obtain the bound
\begin{equation}\label{eta1boundprel}
    \| \beeta_1(e)\|_{\bL^2(\es(e))} \lesssim h_e  \sum_{v\in\Ve} \|\ZZ_p^0(v)\|_{L^2(\es(v))}.
\end{equation}
We extend $\beeta_1(e)$ by zero outside $\es(e)$.

Next, we define $\bpsi_1(e) \in \Zz^{\perp} \bV_p^1(\calTee)$ (a piecewise N\'ed\'elec polynomial of order $p$ on the extended edge star $\es(e)$, orthogonal to all curl-free piecewise N\'ed\'elec polynomials of order $p$ on $\calTee$) such that
\begin{equation}\label{605.1}
\bl \mu_{e} \, \curl \, \bpsi_1(e), \curl \, \bu\br_{\es(e)} =\phi_e(\bu) -\bl \beeta_1(e), \bu \br_{\es(e)}, \quad  \forall \bu \in   \Zz^{\perp} \bV_p^1(\calTee),
\end{equation}
with $\mu_e$ defined in~\eqref{eq:def_mu}.
This problem is well-posed owing to the Poincar\'e inequality~\eqref{onto1}.
Let $\bu \in \Zz \bV_p^1(\calTee)$, so that $\curl \, \bu=\bzero$. Then by the exactness of the discrete sequence~\eqref{localcomplex} from Proposition~\ref{exactcomplex}, we have $\bu=\grad \, m$ for some $m \in V_p^0(\calTee)$. We have
\begin{alignat*}{2}
\bl \beeta_1(e), \grad \, m \br_{\es(e)}&= -\bl \dive \, \beeta_1(e),  m \br_{\es(e)} \qquad && \text{integration by parts, $\beeta_1(e)\in \mbVp{2}(\calTeeA)$ } \\
&=  \sum_{v\in\Ve} \iota_{ev} \bl \ZZ_p^0(v),  m \br_{\es(v)} \qquad && \text{by~\eqref{eta1def}} \\
&= \sum_{v\in\Ve} \iota_{ev} \phi_v(m)   \qquad && \text{by~\eqref{ZZb0}} \\
&= \phi_e (\grad \, m) \qquad && \text{by~\eqref{eq:dofs_iota_0}}. 
\end{alignat*}
Since \eqref{605.1} also holds for all $\bu \in \Zz \bV_p^1(\calTee)$, it holds altogether for all piecewise N\'ed\'elec polynomials of order $p$ on $\es(e)$, i.e., we have
\begin{equation}\label{805.1}
\bl \mu_{e} \, \curl \, \bpsi_1(e), \curl \, \bu\br_{\es(e)}=\phi_e(\bu) -\bl \beeta_1(e), \bu \br_{\es(e)}, \quad  \forall \bu \in \bV_p^1(\calTee).
\end{equation}

Finally, we define 
\begin{equation}\label{defZ1}
\boxed{\bZZ_p^1(e):=\beeta_1(e)+\curl \big(\mu_e \, \curl \, \bpsi_1(e) \big) \quad \text{in $\es(e)$},}
\end{equation}
and extend $\bZZ_p^1(e)$ by zero outside $\es(e)$. 
Notice that each component of $\mu_e \curl \, \bpsi_1(e)$ is a continuous piecewise polynomial of degree $(p+1)$ on the mesh $\ThA$. Indeed, $\bpsi_1(e)$ being a piecewise N\'ed\'elec polynomial of order $p$ on $\es(e)$, $\curl \, \bpsi_1(e)$ is a piecewise polynomial of degree $p$ on the original tetrahedral mesh $\calTee$ (it is a divergence-free Raviart--Thomas of order $p$). Multiplying by $\mu_e$ increases the polynomial degree by one, brings in the Alfeld split $\ThA$, and ensures the continuity of each component. Moreover, $\mu_e \curl \, \bpsi_1(e)$ is zero on $\partial(\es(e))$. Thus, we can take the curl of $\mu_e \curl \, \bpsi_1(e)$ and obtain a piecewise Raviart--Thomas polynomial of order $p$ on the Alfeld split $\calTeeA$ with zero normal component on $\partial\es(e)$. Hence, we indeed have $\bZZ_p^1(e)\in \mbVp{2}(\calTeeA)$.  
It remains to prove~\eqref{ZZc} for $l=1$, \eqref{ZZb1}, and~\eqref{ZZa0}.

(1) Proof of ~\eqref{ZZc} for $l=1$. By construction, $\bZZ_p^1(e)$ 
has support in $\clos(\es(e))$. First, the bound~\eqref{eta1boundprel}, together with~\eqref{ZZc} for $l=0$ and the shape-regularity of the mesh, imply that
\begin{equation}\label{eta1bound}
\| \beeta_1(e)\|_{\bL^2(\es(e))} \lesssim h_e \sum_{v\in\Ve} \| \ZZ_p^0(v)\|_{L^2(\es(v))} 
\lesssim h_e^{-\frac32+1}.
\end{equation}
Using~\eqref{605.1}, we obtain 
\begin{align*}
\| \curl \, \bpsi_1(e)\|_{\bL^2(\es(e))}^2 &\lesssim \bl\mu_e\curl \, \bpsi_1(e),\curl \, \bpsi_1(e)\br_{\es(e)} \\
&=\phi_e(\bpsi_1(e)) -\bl \beeta_1(e), \bpsi_1(e) \br_{\es(e)} \\
&\le h_e \|\bpsi_1(e)\|_{\bL^\infty(\es(e))} + \| \beeta_1(e)\|_{\bL^2(\es(e))} \|\bpsi_1(e)\|_{\bL^2(\es(e))}.
\end{align*}
Proceeding as we did above for $\ZZ_p^0$, cf.~\eqref{eq_inv_ineq}, and in particular invoking the Poincar\'e inequality~\eqref{onto1} gives
\[
    \| \curl \, \bpsi_1(e)\|_{\bL^2(\es(e))} \lesssim h_e^{-\frac32+2}. 
\]
Finally, invoking an inverse estimate, we obtain 
\begin{equation*}
  \| \curl \big(\mu_e \, \curl \, \bpsi_1(e) \big)\| _{\bL^2(\es(e))}\lesssim  h_e^{-\frac32+1}.
\end{equation*}
Combining the above bounds proves~\eqref{ZZc} for $l=1$.

(2) Proof of~\eqref{ZZb1}. 
Using integration by parts shows that, for all $\bu\in \bV_p^1(\calTee)$,
\begin{align*}
\bl \bZZ_p^1(e), \bu \br_{\es(e)} &= \bl \beeta_1(e)+ \curl\big( \mu_e \, \curl \, \bpsi_1(e)\big), \bu \br_{\es(e)} \\
&= \bl \beeta_1 (e), \bu \br_{\es(e)} + \bl  \mu_e \, \curl \, \bpsi_1(e), \curl \,  \bu \br_{\es(e)} = \phi_e(\bu),
\end{align*}
owing to~\eqref{805.1}. Since $\bZZ_p^1(e)$ is zero outside $\es(e)$, the above identity holds for all $\bu \in \bV_p^1(\Th)$. This proves~\eqref{ZZb1}.

(3) Proof of~\eqref{ZZa0}. We readily see that
\begin{equation} \label{eq:div_Z1}
-\dive \, \bZZ_p^1(e)= -\dive \, \beeta_1(e)= \sum_{v\in\Ve} \iota_{ev} \ZZ_p^0(v),
\end{equation}
since $\div\curl=0$ and owing to~\eqref{eta1def}. This proves~\eqref{ZZa0}.

\subsection{Construction and properties of $\bZZ_p^2(f)$}\label{sec:ZZ2}
Let $f \in \mathcal{F}_h$. Owing to~\eqref{eq:div_Z1}, we infer that
\begin{align*}
\dive \left( \sum_{e\in\Ef} \iota_{fe} \bZZ^1_p(e) \right) &=
-\sum_{e\in\Ef} \iota_{fe} \left( \sum_{v\in\Ve} \iota_{ev} \ZZ_p^0(v) \right) \\
&=-\sum_{v\in \Vf} \left(\sum_{e\in\Ev\cap\Ef} \iota_{fe}\iota_{ev}\right) \ZZ_p^0(v) = 0,
\end{align*}
where the last equality follows from~\eqref{eq:CURL_GRAD}. Hence, reasoning as above in~\rev{\eqref{eq_tilde_beeta}} %--\eqref{eq_RT_min} 
and using the Poincar\'e inequality~\eqref{onto1m}, there exists  $\beeta_2(f) \in \Zz^{\perp} \mbVp{1}(\calTefA)$, i.e., a piecewise N\'ed\'elec polynomial of order $p$ on the Alfeld split of the extended face star $\es(f)$, orthogonal to all curl-free piecewise N\'ed\'elec polynomials of order $p$ on $\calTefA$ and with zero tangential component on $\partial \es(f)$, such that 
\begin{equation}\label{eta2def}
 \curl \, \beeta_2(f)= \sum_{e\in\Ef} \iota_{fe}  \bZZ_p^1(e) \text{ on $\es(f)$},
\end{equation}
with the bound
\begin{equation}\label{eta2boundprel}
    \| \beeta_2(f)\|_{\bL^2(\es(f))} \lesssim h_f \sum_{e\in\Ef} \|\bZZ_p^1(e)\|_{\bL^2(\es(e))}.
\end{equation}
We extend $\beeta_2(f)$ by zero outside of $\es(f)$.

Next, we define  $\bpsi_2(f) \in \Zz^{\perp} \bV_p^2(\calTef)$ (a piecewise Raviart--Thomas polynomial of order $p$ on the extended face star $\es(f)$, orthogonal to all divergence-free piecewise Raviart--Thomas polynomials of order $p$ on $\calTef$) such that
\begin{equation}\label{605.2}
\bl \mu_{f} \, \dive \, \bpsi_2(f), \dive \, \bu\br_{\es(f)} =\phi_f(\bu) -\bl \beeta_2(f), \bu \br_{\es(f)}, \quad  \forall \bu \in   \Zz^{\perp} \bV_p^2(\es(f)),
\end{equation}
with $\mu_f$ defined in~\eqref{eq:def_mu}.
This problem is well-posed owing to the Poincar\'e inequality~\eqref{onto2}. 
Let $\bu \in \Zz \bV_p^2(\calTef)$, so that $\dive \, \bu=0$. Then by the exactness of the discrete sequence~\eqref{localcomplex} from Proposition~\ref{exactcomplex}, we have $\bu=\curl \, \bbm$ for some $\bbm\in \bV_p^1(\calTef)$. We have
\begin{alignat*}{2}
\bl \beeta_2(f), \curl \, \bbm \br_{\es(f)}&= \bl \curl \, \beeta_2(f),  \bbm \br_{\es(f)} \qquad && \text{integration by parts, $\beeta_2(f) \in \mbVp{1}(\calTefA)$ } \\
&= \sum_{e\in\Ef} \iota_{fe}  \bl \bZZ_p^1(e), \bbm \br_{\es(f)} \qquad && \text{by~\eqref{eta2def}} \\
&= \sum_{e\in\Ef} \iota_{fe}  \phi_e(\bbm) \qquad && \text{by~\eqref{ZZb1}} \\
&= \phi_f(\curl \, \bbm) \qquad &&\text{by~\eqref{eq:dofs_iota_1}}.
\end{alignat*}
Hence, \eqref{605.2} also holds for all $\bu \in \Zz \bV_p^2(\calTef)$, so that
\begin{equation}\label{805.2}
\bl \mu_f \, \dive \, \bpsi_2(f), \dive \, \bu\br_{\es(f)}=\phi_f(\bu) -\bl \beeta_2(f), \bu \br_{\es(f)}, \quad  \forall \bu \in \bV_p^2(\calTef).
\end{equation}

Finally, we define 
\begin{equation}\label{defZ2}
\boxed{\bZZ_p^2(f):=\beeta_2(f)- \grad \, \big(\mu_f \, \dive \, \bpsi_2(f) \big) \quad\text{in $\es(f)$},}
\end{equation}
and extend $\bZZ_p^2(f)$ by zero outside $\es(f)$. 
Notice that $\mu_f \dive \, \bpsi_2(f)$ is a continuous piecewise polynomial of degree $(p+1)$ on the mesh $\ThA$ and is zero on $\partial(\es(f))$. Thus, we can take the gradient of $\mu_f \dive \, \bpsi_2(f)$ and obtain a piecewise N\'ed\'elec polynomial of order $p$ on the Alfeld split $\calTefA$ with zero tangential component on $\partial\es(f)$. Hence, we indeed have $\bZZ_p^2(f)\in \mbVp{1}(\calTefA)$. 
It remains to prove~\eqref{ZZc} for $l=2$, \eqref{ZZb2}, and~\eqref{ZZa1}.

(1) Proof of~\eqref{ZZc} for $l=2$. By construction, $\bZZ_p^2(f)$ has support in
$\clos(\es(f))$. First, the bound~\eqref{eta2boundprel}, together with~\eqref{ZZc} for $l=1$
and the shape-regularity of the mesh, imply that 
\begin{equation}\label{eta2bound}
\| \beeta_2(f)\|_{\bL^2(\es(f))} \lesssim h_f \sum_{e\in\Ef} \| \bZZ_p^1(e)\|_{\bL^2(\es(e))} \lesssim h_f^{-\frac32+2}.
\end{equation}
Using~\eqref{605.2}, we obtain 
\begin{align*}
\| \dive \, \bpsi_2(f)\|_{L^2(\es(f))}^2 &\lesssim \bl \mu_f\dive \, \bpsi_2(f),\dive \, \bpsi_2(f)\br_{\es(f)} \\
&= \phi_f(\bpsi_2(f)) - \bl \beeta_2(f),\bpsi_2(f)\br_{\es(f)} \\
&\le h_f^2 \|\bpsi_2(f)\|_{\bL^\infty(\es(f))} + \| \beeta_2(f)\|_{\bL^2(\es(f))}\| \bpsi_2(f)\|_{\bL^2(\es(f))}.
\end{align*}
Proceeding as we did above for $\bZZ_p^1(e)$, and in particular invoking the Poincar\'e inequality~\eqref{onto2}, gives
\[
\| \dive \, \bpsi_2(f)\|_{L^2(\es(f))} \lesssim h_f^{-\frac32+3}.
\]
Finally, invoking an inverse estimate, we obtain
\[
\|\grad \big(\mu_f \, \dive \, \bpsi_2(f) \big)\|_{\bL^2(\es(f))} \lesssim h_f^{-\frac32+2}.
\]
Combining the above bounds proves~\eqref{ZZc} for $l=2$.

(2) Proof of~\eqref{ZZb2}. Using integration by parts shows that, for all
$\bu\in\bV_p^2(\calTef)$,
\begin{align*}
\bl \bZZ_p^2(f), \bu \br_{\es(f)} &= \bl \beeta_2(f) - \grad\big( \mu_f \, \dive \, \bpsi_2(f)\big), \bu \br_{\es(f)} \\
& = \bl \beeta_2 (f), \bu \br_{\es(f)}+ \bl  \mu_f \, \dive \, \bpsi_2(f), \dive \, \bu \br_{\es(f)} 
= \phi_f(\bu),
\end{align*}
owing to~\eqref{805.2}. Since $\bZZ_p^2(f)$ is zero outside $\es(f)$, the above identity holds for all $\bu\in\bV_p^2(\Th)$. This proves~\eqref{ZZb2}.

(3) Proof of~\eqref{ZZa1}. We readily see that 
\begin{equation} \label{eq:curl_Z2}
\curl \, \bZZ_p^2(f)= \curl \, \beeta_2(f)= \sum_{e\in\Ef} \iota_{fe} \bZZ_p^1(e),
\end{equation}
since $\curl \,\grad=\bzero$ and owing to~\eqref{eta2def}. This proves~\eqref{ZZa1}. 

\subsection{Construction and properties of $\ZZ_p^3(\tau)$}\label{sec:ZZ3}

Let $\tau\in \Th$. Owing to~\eqref{eq:curl_Z2}, we infer that
\begin{align*}
\curl \left( \sum_{f\in\Ft} \iota_{\tau f} \bZZ_p^2(f)\right) = \sum_{f\in\Ft}\sum_{e\in\Ef}\iota_{\tau f}\iota_{fe}\bZZ_p^1(e) = \sum_{e\in\Et}\sum_{f\in\Fe\cap\Ft} \iota_{\tau f}\iota_{fe}\bZZ_p^1(e) = \bzero,
\end{align*}
where we used~\eqref{eq:DIV_CURL} in the last equality. Hence, by the exactness of the discrete sequence~\eqref{localcomplexm} from Proposition~\ref{exactcomplex} and the Poincar\'e inequality~\eqref{onto0m}, there exists $\eta_3(\tau) \in \Zz^{\perp} \mVp{0}(\calTetA)$, a continuous piecewise polynomial of degree $(p+1)$ on the Alfeld split of the extended element star $\es(\tau)$ with zero value on $\partial \es(\tau)$ (notice that $\Zz^{\perp} \mVp{0}(\calTetA) = \mVp{0}(\calTetA)$), such that 
\begin{equation}\label{eta3def}
-\grad \, \eta_3(\tau)= \sum_{f\in\Ft} \iota_{\tau f} \bZZ_p^2(f) \quad \text{ on $\es(\tau)$},
\end{equation}
with the bound
\begin{equation}\label{eta3boundprel}
\| \eta_3(\tau)\|_{L^2(\es(\tau))} \lesssim h_\tau \sum_{f\in\Ft}  \|\bZZ_p^2(f)\|_{\bL^2(\es(f))}.
\end{equation}
We extend $\eta_3(\tau)$ by zero outside $\es(\tau)$. 

Let $u \in V_p^3(\calTet)$. By the exactness of the discrete sequence~\eqref{localcomplex} from Proposition~\ref{exactcomplex}, there exists $\bbm \in \bV_p^2(\calTet)$ such that $\dive \, \bbm=u$, and we have
\begin{alignat*}{2}
\bl \eta_3(\tau), \dive \,\bbm \br_{\es(\tau)} &=
-\bl \grad \, \eta_3(\tau), \bbm \br_{\es(\tau)} \qquad && \text{integration by parts, $\eta_3(\tau) \in \mVp{0}(\calTetA)$} \\
&= \sum_{f\in\Ft} \iota_{\tau f} \bl \bZZ_p^2(f), \bbm \br_{\es(\tau)}   \qquad && \text{by~\eqref{eta3def}} \\
&= \sum_{f\in\Ft} \iota_{\tau f} \phi_f(\bbm) \qquad && \text{by~\eqref{ZZb2}} \\
&= \phi_\tau(\dive \, \bbm) \qquad && \text{by~\eqref{eq:dofs_iota_2}}.
\end{alignat*}
This proves that
\begin{equation}\label{805.3}
\bl \eta_3(\tau), u \br_{\es(\tau)}= \phi_\tau(u), \quad \forall u \in V_p^3(\calTet).    
\end{equation}

Finally, we simply put
\begin{equation}\label{defZ3}
\boxed{\ZZ_p^3(\tau):=\eta_3(\tau) \quad \text{in $\es(\tau)$},}
\end{equation}
so that $\ZZ_p^3(\tau)\in \mVp{0}(\calTetA)$.
We extend $\ZZ_p^3(\tau)$ by zero outside $\es(\tau)$. 
It remains to prove~\eqref{ZZc} for $l=3$, \eqref{ZZb3}, and~\eqref{ZZa2}. 

(1) By construction, $\ZZ_p^3(\tau)$ has support in $\clos(\es(\tau))$. Moreover, the bound in~\eqref{ZZc} for $l=3$ follows from~\eqref{eta3boundprel} and~\eqref{ZZc} for $l=2$.

(2) \eqref{ZZb3} readily results from~\eqref{805.3}.

(3) \eqref{ZZa2} readily follows from~\eqref{eta3def}.

\section{Local $L^2$-bounded commuting projections}
\label{sec:projections}

In this section, we build the local $L^2$-bounded commuting projections from the local weight
functions devised in Theorem~\ref{th:weights}. The way to do this follows from the ideas in~\cite{Arn_Guz_loc_stab_L2_com_proj_21}. We present here the construction of the projections in the lowest-order and higher-order cases, but for brevity, we only prove their properties in the lowest-order case.

\subsection{Lowest-order projections}
\label{sec:lowest}

In this section, we state and prove our main result on 
the lowest-order projections, thus considering the polynomial degree $p=0$. 
% For all $l\in\{0{:}3\}$, we define the projection $\Pi_0^l:L^2(\Omega;\R^q)\rightarrow V_0^l(\Th)$ with $q=1$ for $l\in\{0,3\}$ and $q=d=3$ for $l\in\{1,2\}$. 
Specifically, we define
$\Pi_0^0:L^2(\Omega)\rightarrow V_0^0(\Th)$,
$\bPi_0^1:\bL^2(\Omega)\rightarrow \bV_0^1(\Th)$,
$\bPi_0^2:\bL^2(\Omega)\rightarrow \bV_0^2(\Th)$,
$\Pi_0^3:L^2(\Omega)\rightarrow V_0^3(\Th)$ as follows:
\begin{subequations} \label{eq:Pi}
\begin{alignat}{2}
\Pi_0^0 (u) &:= \sum_{v \in \Vh} \bl  \ZZ_0^0(v), u \br_{\es(v)} W_v, \qquad &&\forall u\in L^2(\Omega),\label{Pi0} \\ 
\bPi_0^1 (\bu)&:= \sum_{e \in \Eh} \bl  \bZZ_0^1(e), \bu \br_{\es(e)} \bW_e, \qquad &&\forall \bu\in \bL^2(\Omega), \label{Pi1}\\
\bPi_0^2 (\bu)&:= \sum_{f \in \Fh} \bl  \bZZ_0^2(f), \bu \br_{\es(f)} \bW_f, \qquad &&\forall \bu\in \bL^2(\Omega), \label{Pi2}\\ 
\Pi_0^3 (u)&:= \sum_{\tau \in \Th} \bl  \ZZ_0^3(\tau), u \br_{\es(\tau)} W_\tau, \qquad &&\forall u\in L^2(\Omega),   \label{Pi3}
\end{alignat}
\end{subequations}
where the Whitney forms $(W_v)_{v\in\Vh}$, $(\bW_e)_{e\in\Eh}$, $(\bW_f)_{f\in\Fh}$, $(W_\tau)_{\tau\in\Th}$ are defined in Section~\ref{sec:Whitney}.

\begin{theorem}[Lowest-order projections]\label{th:projection}
For all $l\in\{0{:}3\}$, the linear operators $\Pi_0^l$ defined in~\eqref{eq:Pi} satisfy the following properties: \\
\textup{(i)} They are projections, i.e., $\Pi_0^l(u) = u$ for all $u\in V_{0}^l(\Th)$. \\
\textup{(ii)} They are locally $L^2$-bounded: There exists a
constant $C_\Pi>0$, only depending on the mesh shape-regularity parameter $\rho_{\Th}$, such that, for all $l\in\{0{:}3\}$,  
\begin{equation}\label{PIbound}
\| \Pi_0^l (u)\|_{L^2(\tau)} \le C_\Pi \|u\|_{L^2(\es(\tau))}, \quad \forall  \tau \in \Th. \end{equation}
\textup{(iii)} They commute with the differential operators:
\begin{subequations}
\begin{alignat}{2}
\grad (\Pi_0^0 (u))&= \bPi_0^1 (\grad \, u), \qquad && \forall u \in V^0,   \label{commute0}\\
\curl (\bPi_0^1 (\bu))&= \bPi_0^2 (\curl \, \bu), \qquad && \forall \bu \in \bV^1,  \label{commute1}\\
\dive (\bPi_0^2 (\bu))&= \Pi_0^3 (\dive \, \bu),  \qquad && \forall \bu \in \bV^2, \label{commute2}
\end{alignat}
\end{subequations}
with the graph spaces defined in~\eqref{eq:graph_spaces}.
\end{theorem}

\begin{proof}
\rev{The proof follows the arguments of \cite[Theorem~4.3]{Arn_Guz_loc_stab_L2_com_proj_21}.
We only sketch the proof of (iii) since the argument will need to be slightly 
adapted in the more general setting of
Theorem~\ref{th:projectionhom} where a boundary prescription is additionally 
enforced. For brevity, we only detail the proof of~\eqref{commute0}.}
% \textup{(i)} Proof that $\Pi_0^l$ is a projection. We only prove the result for $l=0$; the proofs for 
% $l\ge1$ are similar. Let $u\in V_{0}^0(\Th)$. Then $u=\sum_{v\in\Vh} \phi_v(u)W_v$, and we infer from~\eqref{ZZb0} with $p=0$ that
% \[
% \Pi_0^0(u) = \sum_{v \in \Vh} \bl  \ZZ_0^0(v), u \br_{\es(v)} W_v = \sum_{v\in\Vh} \phi_v(u)W_v = u.
% \]
%
% \textup{(ii)} Proof that $\Pi_0^l$ is locally $L^2$-bounded. Again, we only prove the result for $l=0$; 
% the proofs for $l\ge1$ are similar. Let $\tau\in\Th$. For all $u\in L^2(\Omega)$, we have
% $\Pi_0^0(u)|_\tau = \sum_{v \in \Vt} \bl  \ZZ_0^0(v), u \br_{\es(v)} W_v|_\tau$. Therefore, the
% triangle inequality, the Cauchy--Schwarz inequality, the bound~\eqref{CW} on the Whitney forms, and the bound~\eqref{ZZc} for $l=0$ on the local weight functions imply that
% \begin{align*}
% \|\Pi_0^0(u)\|_{L^2(\tau)} &\le \sum_{v \in \Vt} \|\ZZ_0^0(v)\|_{L^2(\es(v))} \|u\|_{L^2(\es(v))} \|W_v\|_{L^2(\tau)} \\
% &\le \sum_{v \in \Vt} C_Z h_v^{-\frac32+0} C_W h_v^{\frac32-0} \|u\|_{L^2(\es(v))} \lesssim \|u\|_{L^2(\es(\tau))},
% \end{align*}
% where the last bound follows from the mesh shape-regularity.
%
% \textup{(iii)} Proof of commuting property with differential operators. We only prove~\eqref{commute0}, since the proof of the other two identities uses similar arguments. 
Let $u\in V^0$. Since both
$\grad (\Pi_0^0 (u))$ and $\bPi_0^1 (\grad \, u)$ are in $\bV_0^1(\Th)$, it suffices to show that, for all $e\in \Eh$,
\[
\phi_e(\grad (\Pi_0^0 (u))) = \phi_e(\bPi_0^1 (\grad \, u)).
\]  
On the one hand, \eqref{eq:dof_Whitney} and~\eqref{Pi1} show that
\[
\phi_e(\bPi_0^1 (\grad \, u)) = \bl \bZZ_0^1(e), \grad \, u \br_{\es(e)}
= - \bl \dive \, \bZZ_0^1(e), u\br_{\es(e)},
\]
where we used integration by parts and $\bZZ_0^1(e) \in \mbVp{2}(\calTeeA)$. 
On the other hand, using~\eqref{eq:dofs_iota_0}, we have
\[
\phi_e(\grad (\Pi_0^0 (u))) = \sum_{v\in \Ve} \iota_{ev} \phi_v(\Pi_0^0 (u))
= \sum_{v\in \Ve} \iota_{ev} \bl \ZZ_0^0(v), u\br_{\es(v)},
\]
where we used again~\eqref{eq:dof_Whitney} and~\eqref{Pi0} in the last step. We conclude by invoking~\eqref{ZZa0} (for $p=0$)
and observing that $\ZZ_0^0(v)$ is supported in $\clos(\es(v))\subset \clos(\es(e))$ for all $v\in \Ve$.
\end{proof}

\subsection{Higher-order projections}
\label{sec:higher}

We now consider a general polynomial degree $p\ge0$. 
% For all $l\in\{0{:}3\}$, we define the projection $\Pi_p^l:L^2(\Omega;\R^q)\rightarrow V_p^l(\Th)$ with $q=1$ for $l\in\{0,3\}$ and $q=d=3$ for $l\in\{1,2\}$. Specifically, we have
We first define the linear operators
$P_p^l:L^2(\Omega)\rightarrow V_0^l(\Th)$, $l\in\{0{:}3\}$, as follows:
\begin{subequations} \label{eq:Pip}
\begin{alignat}{2}
P_p^0 (u) &:= \sum_{v \in \Vh} \bl  \ZZ_p^0(v), u \br_{\es(v)} W_v, \qquad &&\forall u\in L^2(\Omega),\label{Pip0} \\ 
\bP_p^1 (\bu)&:= \sum_{e \in \Eh} \bl  \bZZ_p^1(e), \bu \br_{\es(e)} \bW_e, \qquad &&\forall \bu\in \bL^2(\Omega), \label{Pip1}\\
\bP_p^2 (\bu)&:= \sum_{f \in \Fh} \bl  \bZZ_p^2(f), \bu \br_{\es(f)} \bW_f, \qquad &&\forall \bu\in \bL^2(\Omega), \label{Pip2}\\ 
P_p^3 (u)&:= \sum_{\tau \in \Th} \bl  \ZZ_p^3(\tau), u \br_{\es(\tau)} W_\tau, \qquad &&\forall u\in L^2(\Omega).   \label{Pip3}
\end{alignat}
\end{subequations}
The only difference with~\eqref{eq:Pi} is that we are now using the higher-order local weight functions. However, the linear operators $P_p^l$ still map onto the lowest-order polynomial spaces $V_0^l(\Th)$. To define the projections
$\Pi_p^0:L^2(\Omega)\rightarrow V_p^0(\Th)$,
$\bPi_p^1:\bL^2(\Omega)\rightarrow \bV_p^1(\Th)$,
$\bPi_p^2:\bL^2(\Omega)\rightarrow \bV_p^2(\Th)$,
$\Pi_p^3:L^2(\Omega)\rightarrow V_p^3(\Th)$, we set, for $l\in\{0,3\}$,
\begin{equation} \label{eq:def_Ppl}
\Pi_p^l(u) := P_p^l(u) + Q_p^l(u-P_p^l(u)), \qquad \forall u\in L^2(\Omega),
\end{equation}
with suitable linear operators $Q_p^l$ defined on $L^2(\Omega)$ and mapping onto
suitable subspaces of $V_p^l(\Th)$. The construction of these operators is described in~\cite{Arn_Guz_loc_stab_L2_com_proj_21}; it hinges on a specific choice of the degrees of freedom in the polynomial
spaces $V_p^l(\Th)$, but does not use the local weight functions specifically. Therefore, we only state our main result concerning the higher-order projections; the proof follows the one
given in~\cite{Arn_Guz_loc_stab_L2_com_proj_21}. 

\begin{theorem}[Higher-order projections]
For all $l\in\{0{:}3\}$, the linear operators $\Pi_p^l$ defined in~\eqref{eq:def_Ppl} satisfy the following properties: \\
\textup{(i)} They are projections, i.e., $\Pi_p^l(u) = u$ for all $u\in V_p^l(\Th)$. \\
\textup{(ii)} They are locally $L^2$-bounded: There exists a
constant $C_\Pi'>0$, only depending on the mesh shape-regularity parameter $\rho_{\Th}$ and the polynomial degree $p$, such that, for all $l\in\{0{:}3\}$,  
\begin{equation}\label{Ppl_bound}
\| \Pi_p^l (u)\|_{L^2(\tau)} \le C_\Pi' \|u\|_{L^2(\es^2(\tau))}, \qquad \forall  \tau \in \Th, \end{equation}
where $\es^2(\tau) := \inte \bigcup_{\tau'\in \Th ; \tau' \cap \es(\tau) \ne \emptyset}\tau'$. \\
\textup{(iii)} They commute with the differential operators:
\begin{subequations}
\begin{alignat}{2}
\grad (\Pi_p^0 (u))&= \bPi_p^1 (\grad \, u), \qquad && \forall u \in V^0,   \label{commute0p}\\
\curl (\bPi_p^1 (\bu))&= \bPi_p^2 (\curl \, \bu), \qquad && \forall \bu \in \bV^1,  \label{commute1p}\\
\dive (\bPi_p^2 (\bu))&= \Pi_p^3 (\dive \, \bu),  \qquad && \forall \bu \in \bV^2. \label{commute2p}
\end{alignat}
\end{subequations}
\end{theorem}

\section{Projections with homogeneous boundary conditions \rev{imposed on part of $\partial\Omega$}}
\label{sec:boundary}

It is desirable to have projections that respect homogeneous boundary conditions \rev{on a part of the boundary $\partial\Omega$}. \rev{Let $\partial\Omega=\Gamma \cup \Gamma^c$ be a partition of the boundary $\partial\Omega$ so that both subsets $\Gamma$ and $\Gamma^c$ are meshed exactly by the underlying simplicial mesh. Both subsets $\Gamma$ and $\Gamma^c$ are supposed to be closed Lipschitz subsets of $\partial\Omega$ and to have disjoint relative interiors. We want to prescribe homogeneous boundary conditions on $\Gamma$. The present construction includes the case where $\Gamma=\partial\Omega$ and $\Gamma^c=\emptyset$. We focus for simplicity on the prescription of homogeneous boundary conditions. We expect that the construction can be further adapted so as to preserve piecewise polynomial data on the boundary, but we postpone this aspect to future work.}

Letting $\bn_\Omega$ denote the unit outward normal to $\Omega$, and recalling the definitions~\eqref{eq:graph_spaces}, the graph spaces with zero boundary conditions \rev{on $\Gamma$} are
\begin{subequations} \begin{align}
  V_{\Gamma}^0 &:= \{ u\in V^0\,:\, u|_{\rev{\Gamma}}=0 \},\\
  \bV_{\Gamma}^1 &:= \{\bu\in \bV^1\,:\, \bu|_{\rev{\Gamma}} {\times} \bn_\Omega=\bzero \}, \\
  \bV_{\Gamma}^2 &:= \{ \bu\in \bV^2\,:\, \bu|_{\rev{\Gamma}} {\cdot}\bn_\Omega =0 \}.
\end{align} \end{subequations}
The corresponding discrete subspaces with homogeneous boundary conditions are  
\begin{equation}
V_{\Gamma,p}^0(\Th):= V_p^0(\Th) \cap V_{\Gamma}^0, \qquad 
\bV_{\Gamma,p}^1(\Th):= \bV_p^1(\Th) \cap \bV_{\Gamma}^1, \qquad
\bV_{\Gamma,p}^2(\Th):= \bV_p^2(\Th) \cap \bV_{\Gamma}^2,
\end{equation}
or, generically, $V_{\Gamma,p}^l(\Th):= V_p^l(\Th) \cap V_{\Gamma}^l$ for all $l\in\{0{:}2\}$.
We also set $V_{\Gamma,p}^3(\Th):=V_p^3(\Th)\cap V_{\Gamma}^3$ with $V_{\Gamma}^3:=\{u\in V^3\,:\,\bl u,1\br_\Omega=0\}$ \rev{when $\Gamma=\partial\Omega$ and $V_{\Gamma}^3=V^3$ otherwise}.  
The projections defined in Section~\ref{sec:projections} do not map $V_{\Gamma}^l$ onto $V_{\Gamma,p}^l(\Th)$. In this section, we define commuting projections that do so and are locally $L^2$-bounded.  

The idea is to modify the weights $\ZZ_p^l(\sigma)$ for some of the geometric entities 
$\sigma \in \Delta_h^l$. Obviously, we want to enforce
$\ZZ_p^l(\sigma)=0$ whenever $\sigma \subset \rev{\Gamma}$, and to keep the weight function $\ZZ_p^l(\sigma)$ unmodified whenever $\partial\es(\sigma)$ does not contain any face located on $\rev{\Gamma}$ (such faces are called \rev{$\Gamma$}-boundary faces). The case where $\sigma\not\subset\rev{\Gamma}$ and $\partial\es(\sigma)$ contains at least one \rev{$\Gamma$}-boundary face is more subtle. For instance, if $e$ is an edge with one \rev{$\Gamma$}-boundary vertex $v_0$ and one interior vertex $v_1$, $\ZZ_p^0(v_0)=0$ has mean zero, whereas $\ZZ_p^0(v_1)$ has mean 1 by~\eqref{eq:int_Z0}. Therefore, the right-hand side of~\eqref{eta1def} cannot have mean zero, and so we do not necessarily have the existence of $\beeta_1(e) \in \Zz^{\perp} \mbVp{2}(\calTeeA) \subset \mbVp{2}(\calTeeA)$  satisfying~\eqref{eta1def}. Thus, the construction presented in Section~\ref{sec:proof} breaks down. The way to fix this is \emph{not to require} that $\beeta_1(e) \in  \mbVp{2}(\calTeeA)$, but that its normal component vanishes only on $\partial\es(e) \setminus \rev{\Gamma}$. 

\subsection{Local polynomial subspaces with boundary prescription}

Let us make the above idea more precise. For all $\sigma \in \Delta_h$, let $\Gamma(\sigma)$
(resp., $\Gac(\sigma)$) 
be the part of the boundary $\partial\es(\sigma)$ which is composed of mesh faces which \rev{are (resp., are not) $\Gamma$-}boundary faces. Notice that $\partial\es(\sigma)=\Gamma(\sigma)\cup \Gac(\sigma)$ and that
$\Gamma(\sigma)\cap \Gac(\sigma)$ has zero \rev{surface} measure.
We partition the geometric entities as
$\Delta_h = \Delta_h\upb\cup \Delta_h\upi \cup \Delta_h\upib$, where
\begin{subequations} \begin{align}
\Delta_h\upb &:= \{\sigma \in \Delta_h\,:\, \sigma \subset \rev{\Gamma}\} ,\\
\Delta_h\upi &:= \{\sigma \in \Delta_h\setminus \Delta_h\upb\,:\, |\Gamma(\sigma)|=0\},\\
\Delta_h\upib &:=\{\sigma \in \Delta_h\setminus \Delta_h\upb\,:\, |\Gamma(\sigma)|>0\},
\end{align} \end{subequations}
and where $|\SCAL|$ denotes the \rev{surface} measure. Specifically, this leads to the following partitions:
\begin{subequations} \begin{align}
\Vh &= \Vh\upb \cup \Vh\upi \cup \Vh\upib, \\
\Eh &= \Eh\upb \cup \Eh\upi \cup \Eh\upib, \\
\Fh &= \Fh\upb \cup \Fh\upi \cup \Fh\upib, \\
\Th &= \phantom{\Th\upb \cup {}} \Th\upi \cup \Th\upib,
\end{align} \end{subequations}
since $\Th\upb = \emptyset$. We assume that $|\Gac(\tau)|\neq 0$ for all $\tau\in\Th$. This is a mild assumption, since $|\Gac(\tau)|= 0$ if and only if $\Gamma(\tau) = \partial\es(\tau)$, which corresponds to the case where $\partial\es(\tau) \subset \Gamma$, which can only occur on very coarse meshes. 

We define the following local polynomial spaces for all $\sigma \in \Delta_h$,
\begin{subequations} \label{eq:def_loc_pol_Gamma} \begin{align}
V_{\Gamma, p}^{0}(\calTes) &:= \{v \in V_p^{0}(\calTes)\,:\, v|_{\Gamma(\sigma)}=0\}, \\
\bV_{\Gamma, p}^{1}(\calTes) &:= \{\bv \in \bV_p^{1}(\calTes)\,:\, \bv|_{\Gamma(\sigma)}\times \bn_{\es(\sigma)}=\bzero\}, \\
\bV_{\Gamma, p}^{2}(\calTes) &:= \{\bv \in \bV_p^{2}(\calTes)\,:\, \bv|_{\Gamma(\sigma)}\SCAL\bn_{\es(\sigma)}=0\},
\end{align} \end{subequations}
as well as 
\begin{subequations} \label{eq:def_loc_pol_Gamma_A} \begin{align}
V_{\Gac, p}^{0}(\calTesA) &:= \{v \in V_p^{0}(\calTesA)\,:\, v|_{\Gac(\sigma)}=0\}, \\
\bV_{\Gac, p}^{1}(\calTesA) &:= \{\bv \in \bV_p^{1}(\calTesA)\,:\, \bv|_{\Gac(\sigma)}\times \bn_{\es(\sigma)}=\bzero\}, \\
\bV_{\Gac, p}^{2}(\calTesA) &:= \{\bv \in \bV_p^{2}(\calTesA)\,:\, \bv|_{\Gac(\sigma)}\SCAL \bn_{\es(\sigma)}=0\},
\end{align} \end{subequations}
where $\bn_{\es(\sigma)}$ denotes the unit outward normal to $\es(\sigma)$.
Notice that $V_{\Gamma, p}^{l}(\calTes) = V_{p}^{l}(\calTes)$ and $V_{\Gac, p}^{l}(\calTesA)
= \mV_{p}^{l}(\calTesA)$ for all $\sigma\in\Delta_h\upi$.

\rev{We assume as above that $\clos(\es(\sigma))$ is contractible and additionally 
that both $\Gamma(\sigma)$ and $\Gac(\sigma)$
have only one simply connected component (i.e., all the boundary faces composing either  
$\Gamma(\sigma)$ or $\Gac(\sigma)$ are edge-connected). This latter assumption 
is satisfied if the mesh is fine enough. Then, invoking 
\cite[Section 3.2.2]{Hipt_Pech_disc_1fo_19} or \cite[Theorem 18]{Hiptmair_99},
we infer that the following sequences are exact:}
\begin{subequations}
\begin{equation}\label{localcomplexPartialBCB}
\mathbb{R}
\stackrel{\subset}{\xrightarrow{\hspace*{0.5cm}}}\
V_{\Gamma, p}^{0}(\calTes)
\stackrel{\grad}{\xrightarrow{\hspace*{0.5cm}}}\
\bV_{\Gamma, p}^{1} (\calTes)
\stackrel{\curl}{\xrightarrow{\hspace*{0.5cm}}}\
\bV_{\Gamma, p}^{2} (\calTes)
\stackrel{\dive}{\xrightarrow{\hspace*{0.5cm}}}\
V_{p}^{3}(\calTes) 
\stackrel{}{\xrightarrow{\hspace*{0.5cm}}}\
0
\end{equation}
\begin{equation}\label{localcomplexPartialBCA}
\mathbb{R}
\stackrel{\subset}{\xrightarrow{\hspace*{0.5cm}}}\
V_{\Gac, p}^{0}(\calTesA)
\stackrel{\grad}{\xrightarrow{\hspace*{0.5cm}}}\
\bV_{\Gac, p}^{1} (\calTesA)
\stackrel{\curl}{\xrightarrow{\hspace*{0.5cm}}}\
\bV_{\Gac, p}^{2} (\calTesA)
\stackrel{\dive}{\xrightarrow{\hspace*{0.5cm}}}\
V_{\Gac, p}^{3}(\calTesA) 
\stackrel{}{\xrightarrow{\hspace*{0.5cm}}}\
0
\end{equation}
\end{subequations}
where $V_{\Gac, p}^{3}(\calTesA):=\{u\in V_p^3(\calTesA)\,:\,\bl u,1\br_{\es(\sigma)}=0\}$ if $\Gac(\sigma)=\partial\es(\sigma)$ and $V_{\Gac, p}^{3}(\calTesA):=V_p^3(\calTesA)$ otherwise.

Finally, we can establish the analogues of the 
discrete Poincar\'e inequalities in Propositions~\ref{prop:discP_star} and~\ref{prop:discP_star_Alfeld} after defining the kernels and orthogonal subspaces of the local polynomial spaces defined
in~\eqref{eq:def_loc_pol_Gamma} and~\eqref{eq:def_loc_pol_Gamma_A}.

\subsection{Main result on weight functions with boundary prescription}

\begin{theorem}[Local weight functions with boundary prescription] \label{th:weightshom}
Let $p\ge0$ be the polynomial degree. 
There exist local weight functions
%a family of local weight functions $\{\ZZ_p^l(\sigma)\}_{\sigma \in \Delta_h^l}$ with $\ZZ_p^l(\sigma) \in \mVp{3-l}(\es(\sigma))$, that is,
$\ZZ_{\Gac,p}^0(v) \in V_p^3(\calTevA)$ for all vertices $v \in \Vh$,
$\bZZ_{\Gac,p}^1(e) \in \bV_{\Gac,p}^{2}(\calTeeA)$ for all edges $e\in \Eh$,
$\bZZ_{\Gac,p}^2(f) \in \bV_{\Gac,p}^{1}(\calTefA)$ for all faces $f\in \Fh$, and
$\ZZ_{\Gac,p}^3(\tau) \in V_{\Gac,p}^{0}(\calTetA)$ for all elements $\tau \in\Th$,
so that $\ZZ_{\Gac,p}^l(\sigma)=\ZZ_p^l(\sigma)$ for all $l\in\{0{:}3\}$ and all $\sigma\in \Delta_h^l\cap \Delta_h\upi$, as well as:\\
\textup{(i)} Nullity for boundary $l$-simplices: For all $l\in\{0{:}3\}$,
\begin{equation}\label{zeroboundarycon}
\ZZ_{\Gac,p}^l(\sigma)=0, \; \forall \sigma\in \Delta_h^l\cap \Delta_h\upb.
\end{equation}
\textup{(ii)} Support and $L^2$-norm: For all $l\in\{0{:}3\}$ and all $\sig \in \Delta_h^l$,
\begin{equation}\label{ZZchom}
\supp\ZZ_{\Gac,p}^l(\sig) \subseteq \clos(\es(\sig)),  \qquad
\|\ZZ_{\Gac,p}^l(\sig)\|_{L^2(\es(\sig))}  \le C_{\mathring{Z}} h_{\sig} ^{-\frac{3}{2}+l}, 
\end{equation}
where $C_{\mathring{Z}}$ only depends on the mesh shape-regularity parameter $\rho_{\Th}$ and the polynomial degree $p$.
\\
\textup{(iii)} Relation to canonical degrees of freedom:
\begin{subequations}\label{ZZbhom}
\begin{alignat}{3}
 \bl \ZZ_{\Gac,p}^0(v), u \br_{\es(v)}&=  \phi_v(u), \qquad && \forall u \in V_{\Gamma,p}^0(\Th), &&\quad \forall v \in  \Delta_h^0=\Vh, \label{ZZb0hom}\\
 \bl \bZZ_{\Gac,p}^1(e), \bu \br_{\es(e)}&=  \phi_e(\bu), \qquad && \forall \bu \in \bV_{\Gamma,p}^1(\Th), &&\quad \forall e \in  \Delta_h^1=\Eh, \label{ZZb1hom} \\
  \bl \bZZ_{\Gac,p}^2(f), \bu \br_{\es(f)}&=  \phi_f(\bu), \qquad && \forall \bu \in \bV_{\Gamma,p}^2(\Th), &&\quad \forall f \in  \Delta_h^2=\Fh,  \label{ZZb2hom}\\
 \bl \ZZ_{\Gac,p}^3(\tau), u \br_{\es(\tau)}&=  \phi_\tau(u),  \qquad && \forall u \in V_{p}^3(\Th), &&\quad \forall \tau \in  \Delta_h^3=\Th. \label{ZZb3hom}
\end{alignat}
\end{subequations}
\textup{(iv)} Relation to differential operators:
\begin{subequations}\label{ZZahom}
\begin{alignat}{2}
-\dive \, \bZZ_{\Gac,p}^1(e)&= \sum_{v\in \Ve} \iota_{ev} \ZZ_{\Gac,p}^0(v), \qquad && \forall e \in \Delta_h^1=\Eh,  \label{ZZa0hom}\\
\curl \, \bZZ_{\Gac,p}^2(f)&= \sum_{e\in \Ef} \iota_{fe} \bZZ_{\Gac,p}^1(e), \qquad && \forall f \in \Delta_h^2=\Fh, \label{ZZa1hom} \\
-\grad \, \ZZ_{\Gac,p}^3(\tau)&= \sum_{f\in \Ft} \iota_{\tau f} \bZZ_{\Gac,p}^2(f), \qquad && \forall \tau \in \Delta_h^3=\Th. \label{ZZa2hom}
\end{alignat} \end{subequations}
\end{theorem}

The proof of Theorem~\ref{th:weightshom} is postponed to Section~\ref{sec:proof_weight_hom}. 

\subsection{Local $L^2$-bounded commuting projections with boundary prescription} 

In this section, we sketch the construction of the lowest-order and higher-order projections. 
For all $l\in\{0{:}3\}$, the lowest-order projections with boundary prescription
$\Pi_{\Gamma,0}^l:L^2(\Omega)\rightarrow V_{\Gamma,0}^l(\Th)$ are defined 
using the formulas~\eqref{eq:Pi}, but employing the weights from Theorem~\ref{th:weightshom} 
with $p=0$.

\begin{theorem}[Lowest-order projections with boundary prescription]\label{th:projectionhom}
For all $l\in\{0{:}3\}$, the linear operators $\Pi_{\Gamma,0}^l$ defined above satisfy the following properties: \\
\textup{(i)} They map $V_{\Gamma}^l$ onto $V_{\Gamma,0}^l(\Th)$ for all $l\in\{0{:}3\}$. \\
\textup{(ii)} They are projections, i.e., $\Pi_{\Gamma,0}^l(u) = u$ for all $u\in V_{\Gamma,0}^l(\Th)$. \\
\textup{(ii)} They are locally $L^2$-bounded: There exists a
constant $C_{\Pi,\Gamma}>0$, only depending on the mesh shape-regularity parameter $\rho_{\Th}$, such that, for all $l\in\{0{:}3\}$,  
\begin{equation}\label{mPIbound}
\| \Pi_{\Gamma,0}^l (u)\|_{L^2(\tau)} \le C_{\Pi,\Gamma} \|u\|_{L^2(\es(\tau))}, \quad \forall  \tau \in \Th. \end{equation}
\textup{(iv)} They commute with the differential operators:
\begin{subequations}
\begin{alignat}{2}
\grad (\Pi_{\Gamma,0}^0 (u))&= \bPi_{\Gamma,0}^1 (\grad \, u), \qquad && \forall u \in V_{\Gamma}^0,   \label{commute0hom}\\
\curl (\bPi_{\Gamma,0}^1 (\bu))&= \bPi_{\Gamma,0}^2 (\curl \, \bu), \qquad && \forall \bu \in \bV_{\Gamma}^1,  \label{commute1hom}\\
\dive (\bPi_{\Gamma,0}^2 (\bu))&= \Pi_{\Gamma,0}^3 (\dive \, \bu),  \qquad && \forall \bu \in \bV_{\Gamma}^2. \label{commute2hom}
\end{alignat}
\end{subequations}
\end{theorem}

\begin{proof} 
We only highlight the two differences in the proof with respect to the proof of Theorem~\ref{th:projection}. 

(1) The fact that (i) holds follows from~\eqref{zeroboundarycon} and the fact that, for a geometric entity $\sigma \in \Delta_h^l$, $l\in\{0{:}3\}$, the corresponding Whitney form is in $V_{\Gamma,0}^l(\Th)$ whenever $\sigma \not\in \Delta_h\upb$.

(2) Consider the commuting property~\eqref{commute0hom}. 
Let $u\in V_{\Gamma}^0$. Since both $\grad (\Pi_{\Gamma,0}^0 (u))$ and $\bPi_{\Gamma,0}^1 (\grad \, u)$ are in $\bV_{\Gamma,0}^1(\Th)$, it suffices to show that, for all $e\in \Eh\setminus \Eh\upb$,
\[
\phi_e(\grad (\Pi_{\Gamma,0}^0 (u))) = \phi_e(\bPi_{\Gamma,0}^1 (\grad \, u)).
\] 
On the one hand, invoking~\eqref{eq:dof_Whitney} and~\eqref{Pi1} (employing $\bZZ_{\Gac,0}^1(e)$ in place of $\bZZ_0^1(e)$) and using integration by parts gives
\[
\phi_e(\bPi_{\Gamma,0}^1 (\grad \, u)) = \bl \bZZ_{\Gac,0}^1(e), \grad \, u \br_{\es(e)}
= - \bl \dive \, \bZZ_{\Gac,0}^1(e), u\br_{\es(e)},
\]
since $u$ vanishes on $\Gamma(e)$ and the normal component of $\bZZ_{\Gac,0}^1(e)$ vanishes on $\Gac(e)$. On the other hand, using~\eqref{eq:dofs_iota_0} followed by~\eqref{eq:dof_Whitney} and~\eqref{Pi0} (employing $\ZZ_{\Gac,0}^0(v)$ in place of $\ZZ_0^0(v)$), we have
\[
\phi_e(\grad (\Pi_{\Gamma,0}^0 (u))) = \sum_{v\in \Ve} \iota_{ev} \phi_v(\Pi_{\Gamma,0}^0 (u))
= \sum_{v\in \Ve} \iota_{ev} \bl \ZZ_{\Gac,0}^0(v), u\br_{\es(v)}.
\]
The commuting property~\eqref{commute0hom} follows from~\eqref{ZZa0hom} (for $p=0$). The proof of the other two commuting properties is similar.
\end{proof}

For all $l\in\{0{:}3\}$, we now define the higher-order projection
$\Pi_{\Gamma,p}^l:L^2(\Omega)\rightarrow V_{\Gamma,p}^l(\Th)$ so that
\begin{equation} \label{eq:def_mPpl}
\Pi_{\Gamma,p}^l(u) := P_{\Gamma,p}^l(u) + Q_{\Gamma,p}^l(u-P_{\Gamma,p}^l(u)), \qquad \forall u\in L^2(\Omega),
\end{equation}
with suitable operators $P_{\Gamma,p}^l$ and $Q_{\Gamma,p}^l$. 
In particular, we can define $P_{\Gamma,p}^l$ using~\eqref{eq:Pip}, but now using the higher-order weight functions $\ZZ_{\Gac,p}^l$. It remains to define the operator $Q_{\Gamma,p}^l$. To this end, we notice that the projections $Q_p^l$ are defined in \cite[Section~5.4]{Arn_Guz_loc_stab_L2_com_proj_21} as a sum over all $\sigma\in\Delta_h^l$ using their own set of weights, say $\mathsf{U}_p^l(\sigma)$, which satisfy analogous properties to the weights $\ZZ_p^l(\sigma)$. However, the weights $\mathsf{U}_p^l(\sigma)$ are defined only in terms of their associated simplex (as opposed to the weights $\ZZ_p^l(\sigma)$ which require the consideration of $\ZZ_p^{l-1}$ on $\es(\sigma)$). Hence, we can define $Q_{\Gamma,p}^l$ following exactly the definition of $Q_p^l$ in \cite[Section~5.4]{Arn_Guz_loc_stab_L2_com_proj_21}, but now summing only over $\sigma\in\Delta_h^l\setminus\Delta_h\upb$. The properties of the resulting projection are stated in the following result. The proof follows again the lines of that given in~\cite{Arn_Guz_loc_stab_L2_com_proj_21} and is skipped for brevity. 

\begin{theorem}[Higher-order projections with boundary prescription]
For all $l\in\{0{:}3\}$, the linear operators $\Pi_{\Gamma,p}^l$ defined above satisfy the following properties: \\
\textup{(i)} They map $V_{\Gamma}^l$ onto $V_{\Gamma,p}^l(\Th)$ for all $l\in\{0{:}3\}$. \\ 
\textup{(ii)} They are projections, i.e., $\Pi_{\Gamma,p}^l(u) = u$ for all $u\in V_{\Gamma,p}^l(\Th)$. \\
\textup{(iii)} They are locally $L^2$-bounded: There exists a
constant $C_{\Pi,\Gamma}'>0$, only depending on the mesh shape-regularity parameter $\rho_{\Th}$ and the polynomial degree $p$, such that, for all $l\in\{0{:}3\}$,  
\begin{equation}\label{mPpl_bound}
\| \Pi_{\Gamma,p}^l (u)\|_{L^2(\tau)} \le C_{\Pi,\Gamma}' \|u\|_{L^2(\es^2(\tau))}, \quad \forall  \tau \in \Th. 
\end{equation}
\textup{(iv)} They commute with the differential operators:
\begin{subequations}
\begin{alignat}{2}
\grad (\Pi_{\Gamma,p}^0 (u))&= \bPi_{\Gamma,p}^1 (\grad \, u), \qquad && \forall u \in V_{\Gamma}^0,   \label{commute0mp}\\
\curl (\bPi_{\Gamma,p}^1 (\bu))&= \bPi_{\Gamma,p}^2 (\curl \, \bu), \qquad && \forall \bu \in \bV_{\Gamma}^1,  \label{commute1mp}\\
\dive (\bPi_{\Gamma,p}^2 (\bu))&= \Pi_{\Gamma,p}^3 (\dive \, \bu),  \qquad && \forall \bu \in \bV_{\Gamma}^2. \label{commute2mp}
\end{alignat}
\end{subequations}
\end{theorem}

\subsection{Proof of Theorem~\ref{th:weightshom}}
\label{sec:proof_weight_hom}

We only highlight the differences with respect to the proof of Theorem~\ref{th:weights} in Section~\ref{sec:proof}.
Recall that, for all $l\in\{0{:}3\}$, we set $\ZZ_{\Gac,p}^l(\sigma)=0$ for all $\sigma \in
\Delta_h^l\cap \Delta_h\upb$, and $\ZZ_{\Gac,p}^l(\sigma)=\ZZ_p^l(\sigma)$ for all $\sigma \in
\Delta_h^l\cap \Delta_h\upi$. Hence, property (i) in Theorem~\ref{th:weightshom}, as well
as properties (ii) and (iii) for all $\sigma\in \Delta_h^l\cap (\Delta_h\upb\cup\Delta_h\upi)$,
are automatically satisfied. It remains to construct the weight functions $\ZZ_{\Gac,p}^l(\sigma)$
for all $\sigma\in \Delta_h^l\cap\Delta_h\upib$ so that they satisfy properties (ii) and (iii), and to verify property (iv). In all cases, property (ii) (support and $L^2$-norm) is proved as in Section~\ref{sec:proof}.

\subsubsection{Construction and properties of $\ZZ_{\Gac,p}^0(v)$}

Let $v\in\Vh\upib$. We define $\mpsi_0(v) \in V_{\Gamma,p}^0(\calTev)=\Zz^\perp V_{\Gamma,p}^0(\calTev)$ by
\begin{equation*}
\bl \mu_{v} \, \grad \, \mpsi_0(v), \grad \, u\br_{\es(v)} = \phi_v(u) -\bl \meta_0(v), u \br_{\es(v)}, \quad  \forall u \in  V_{\Gamma,p}^0(\calTev),
\end{equation*}
with $\meta_0(v):=\eta_0(v)$ (as defined above in~\eqref{eq_eta_0}) and set 
\begin{equation*}
\boxed{\ZZ_{\Gac,p}^0(v):=\meta_0(v)- \dive \big(\mu_v \, \grad \, \mpsi_0(v) \big) \quad \text{in $\es(v)$}.}
\end{equation*}
Notice that $\ZZ_{\Gac,p}^0(v) \in V_p^3(\calTevA)$, but we no longer have $\int_{\es(v)} \ZZ_{\Gac,p}^0(v)=1$ since the constant function $u=1$ is no longer in $V_{\Gamma,p}^0(\calTev)$.
It remains to prove~\eqref{ZZb0hom}. For all $u\in V_{\Gamma,p}^0(\calTev)$, we have
\begin{align*}
\bl \ZZ_{\Gac,p}^0(v), u \br_{\es(v)}&= \bl \meta_0(v)- \dive\big( \mu_v \, \grad \, \mpsi_0(v)\big), u \br_{\es(v)} \\
&= \bl \meta_0 (v), u \br_{\es(v)}+ \bl  \mu_v \, \grad \, \mpsi_0(v), \grad \,  u \br_{\es(v)} = \phi_v(u),
\end{align*}
where we used integration by parts (recall that $\mu_v$ vanishes on $\partial \es(v)$) and the definition of $\mpsi_0(v)$. Since $\ZZ_{\Gac,p}^0(v)$ is zero outside $\es(v)$, this identity holds for all $u\in V_{\Gamma,p}^0(\Th)$.

\subsubsection{Construction and properties of $\bZZ_{\Gac,p}^1(e)$}

Let $e\in\Eh\upib$. Reasoning as in~\eqref{eq_tilde_beeta}, there exists a suitably bounded field $\mbeeta_1(e) \in \Zz^{\perp} \bV_{\Gac,p}^2(\calTeeA)$ such that 
\begin{equation*}
-\dive \, \mbeeta_1(e)= \sum_{v\in\Ve} \iota_{ev} \ZZ_{\Gac,p}^0(v) \quad \text{ on $\es(e)$}.
\end{equation*}
Next, we define $\mbpsi_1(e) \in \Zz^{\perp} \bV_{\Gamma,p}^1(\calTee)$ such that
\begin{equation*}
\bl \mu_{e} \, \curl \, \mbpsi_1(e), \curl \, \bu\br_{\es(e)} =\phi_e(\bu) -\bl \mbeeta_1(e), \bu \br_{\es(e)}, \quad  \forall \bu \in   \Zz^{\perp} \bV_{\Gamma,p}^1(\calTee).
\end{equation*}
Finally, we set 
\begin{equation*}
\boxed{\bZZ_{\Gac,p}^1(e):=\mbeeta_1(e)+\curl \big(\mu_e \, \curl \, \mbpsi_1(e) \big) \quad \text{in $\es(e)$}.}
\end{equation*}
It remains to prove~\eqref{ZZb1hom} and~\eqref{ZZa0hom}.

(1) Proof of~\eqref{ZZb1hom}. The exactness of the sequence~\eqref{localcomplexPartialBCB} implies that any 
$\bu \in \Zz \bV_{\Gamma,p}^1(\calTee)$ is such that $\bu=\grad \, m$ for some $m \in V_{\Gamma,p}^0(\calTee)$. Integration by parts gives
\[
\bl \mbeeta_1(e), \grad \, m \br_{\es(e)}= -\bl \dive \, \mbeeta_1(e),  m \br_{\es(e)},
\]
since $m$ vanishes on $\Gamma(e)$ and the normal component of $\mbeeta_1(e)$ vanishes on
$\Gac(e)$. Proceeding as in Section~\ref{sec:ZZ1} then shows that 
\begin{equation*}
\bl \mu_{e} \, \curl \, \mbpsi_1(e), \curl \, \bu\br_{\es(e)}=\phi_e(\bu) -\bl \mbeeta_1(e), \bu \br_{\es(e)}, \quad  \forall \bu \in \bV_{\Gamma,p}^1(\calTee).
\end{equation*}
Integration by parts (recall that $\mu_e$ vanishes on $\partial\es(e)$) then shows that 
$\bl \bZZ_{\Gac,p}^1(e),\bu\br_{\es(e)} = \phi_e(\bu)$ for all $\bu \in \bV_{\Gamma,p}^1(\calTee)$,
and since $\bZZ_{\Gac,p}^1(e)$ is extended by zero outside $\es(e)$, this proves~\eqref{ZZb1hom}.

(2) \eqref{ZZa0hom} readily follows from $- \dive \, \bZZ_{\Gac,p}^1(e) = -\dive \, \mbeeta_1(e)= \sum_{v\in\Ve} \iota_{ev} \ZZ_{\Gac,p}^0(v)$.

\subsubsection{Construction and properties of $\bZZ_{\Gac,p}^2(f)$}

Let $f\in\Fh\upib$. Reasoning as in~\eqref{eta2def}, there exists a suitably bounded field $\mbeeta_2(f) \in \Zz^{\perp} \bV_{\Gac,p}^1(\calTefA)$ such that 
\begin{equation*}
 \curl \, \mbeeta_2(f)= \sum_{e\in\Ef} \iota_{fe}  \bZZ_{\Gac,p}^1(e) \text{ on $\es(f)$},
\end{equation*}
observing that the right-hand side sits in $\bV_{\Gac,p}^2(\calTefA)$ and is divergence-free
(by the same arguments as in Section~\ref{sec:ZZ2}). 
Next, we define  $\mbpsi_2(f) \in \Zz^{\perp} \bV_{\Gamma,p}^2(\calTef)$ such that
\begin{equation*}
\bl \mu_{f} \, \dive \, \mbpsi_2(f), \dive \, \bu\br_{\es(f)} =\phi_f(\bu) -\bl \mbeeta_2(f), \bu \br_{\es(f)}, \quad  \forall \bu \in   \Zz^{\perp} \bV_{\Gamma,p}^2(\es(f)).
\end{equation*}
Finally, we set
\begin{equation*}
\boxed{\bZZ_{\Gac,p}^2(f):=\mbeeta_2(f)- \grad \big(\mu_f \, \dive \, \mbpsi_2(f) \big) \quad\text{in $\es(f)$}.}
\end{equation*}
It remains to prove~\eqref{ZZb2hom} and~\eqref{ZZa1hom}.

(1) Proof of~\eqref{ZZb2hom}. The exactness of the sequence~\eqref{localcomplexPartialBCB} implies that any $\bu \in \Zz \bV_{\Gamma,p}^2(\calTef)$ is such that $\bu=\curl \, \bbm$ for some 
$\bbm\in \bV_{\Gamma,p}^1(\calTef)$. Integration by parts gives
\[
\bl \mbeeta_2(f), \curl \, \bbm \br_{\es(f)}= \bl \curl \, \mbeeta_2(f),  \bbm \br_{\es(f)},
\]
since the tangential component of $\bbm$ vanishes on $\Gamma(f)$ and the tangential component
of $\mbeeta_2(f)$ vanishes on $\Gac(f)$. Proceeding as in Section~\ref{sec:ZZ2} then shows that 
\[
\bl \mu_{f} \, \dive \, \mbpsi_2(f), \dive \, \bu\br_{\es(f)} =\phi_f(\bu) -\bl \mbeeta_2(f), \bu \br_{\es(f)}, \quad  \forall \bu \in \bV_{\Gamma,p}^2(\es(f)).
\]
Integration by parts (recall that $\mu_f$ vanishes on $\partial\es(f)$) then shows that
$\bl \bZZ_{\Gac,p}^2(f),\bu\br_{\es(f)}=\phi_f(\bu)$ for all $\bu \in \bV_{\Gamma,p}^2(\es(f))$, and 
since $\bZZ_{\Gac,p}^2(f)$ is extended by zero outside $\es(f)$, this proves~\eqref{ZZb2hom}.

(2) \eqref{ZZa1hom} readily follows from 
$\curl \, \bZZ_{\Gac,p}^2(f) = \curl \, \mbeeta_2(f)= \sum_{e\in\Ef} \iota_{fe} \bZZ_{\Gac,p}^1(e)$.

\subsubsection{Construction and properties of $\ZZ_{\Gac,p}^3(\tau)$}

Let $\tau\in\Th\upib$. The exactness of the sequence~\eqref{localcomplexPartialBCA} implies the existence of $\meta_3(\tau) \in \Zz^{\perp} V_{\Gac,p}^0(\calTetA)$ such that 
\begin{equation*}
-\grad \, \meta_3(\tau)= \sum_{f\in\Ft} \iota_{\tau f} \bZZ_{\Gac,p}^2(f) \quad \text{ on $\es(\tau)$},
\end{equation*}
observing that the right-hand side sits in $\bV_{\Gac,p}^1(\calTetA)$ and is curl-free (by the same arguments as in Section~\ref{sec:ZZ3}). Again, $\meta_3(\tau)$ is suitably bounded owing to the Poincar\'e inequality. Finally, we set
\begin{equation*}
\boxed{\ZZ_{\Gac,p}^3(\tau):=\meta_3(\tau) \quad \text{in $\es(\tau)$}.}
\end{equation*}
It remains to prove~\eqref{ZZb3hom} and~\eqref{ZZa2hom}.

(1) Proof of~\eqref{ZZb3hom}. The exactness of the sequence~\eqref{localcomplexPartialBCB} implies that any $u\in V_{p}^3(\calTet)$ is such that $u=\dive \, \bbm$ for some $\bbm\in \bV_{\Gamma,p}^2(\calTet)$. Notice, in particular, that $|\Gac(\tau)|\neq 0$ by assumption. Integration by parts gives
\[
\bl \meta_3(\tau), \dive \,\bbm \br_{\es(\tau)} = -\bl \grad \, \meta_3(\tau), \bbm \br_{\es(\tau)},
\]
since the normal component of $\bbm$ vanishes on $\Gamma(\tau)$ and $\meta_3(\tau)$ vanishes on $\Gac(\tau)$. Proceeding as in Section~\ref{sec:ZZ3} then shows that
\[
\bl \meta_3(\tau),u\br_{\es(\tau)} = \phi_\tau(u), \qquad \forall u\in V_{p}^3(\calTet).
\]
This proves~\eqref{ZZb3hom}.

(2) \eqref{ZZa2hom} readily follows from
$- \grad \, \ZZ_{\Gac,p}^3(\tau) = - \grad \, \meta_3(\tau) = \sum_{f\in\Ft} \iota_{\tau f} \bZZ_{\Gac,p}^2(f)$.

\appendix

\section{Proof of discrete Poincar\'e inequalities (Propositions~\ref{prop:discP_star} and~\ref{prop:discP_star_Alfeld})}
\label{sec:proof_Poinc}

In this appendix, we prove Propositions~\ref{prop:discP_star} and~\ref{prop:discP_star_Alfeld}. 
Since the kernel of the gradient operator is trivial, the discrete Poincar\'e inequalities~\eqref{onto0} and~\eqref{onto0m} are straightforward consequences of their respective
continuous counterparts.
Thus, it only remains to prove the discrete Poincar\'e inequalities involving the curl and divergence operators, i.e., \eqref{onto1}--\eqref{onto2} and~\eqref{onto1m}--\eqref{onto2m}.
We present a proof relying on a constructive argument.
Other arguments can be invoked to prove the discrete Poincar\'e inequalities, as discussed in~\cite{PGEV25}. \rev{The present proof is different from the one given in \cite[Sec.~6.3]{PGEV25}; it shares the idea of considering a finite set of reference patches and piecewise Piola maps on each reference patch, but employs a different argument to conclude.}

\subsection{Unified formalism}
To avoid the proliferation of cases, we adopt a unified formalism. For all $\sigma \in \Delta_h$, we let $\omega:=\es(\sigma)$, and use $\| \SCAL \|_{L^2(\omega)}$ to
generically refer to the $L^2(\omega)$-norm of functions or fields depending on the context. 
Notice that $h_\omega \lesssim h_\sigma$.
In addition, with a slight abuse of notation, we let $\Tom$ denote either $\calTes$ or $\calTesA$ depending on context. 
We set $d^1:=\curl$, $d^2:=\dive$, 
and depending on whether boundary conditions are enforced or not, we set, 
for all $l\in\{1{:}2\}$, ${\widetilde V}_p^l(\Tom):=V_p^l(\calTes)$ or $\mV_p^l(\calTesA)$ and
\begin{subequations} \begin{align}
\Zz {\widetilde V}_p^l(\Tom) &:= \{ u \in {\widetilde V}_p^l(\Tom) : d^lu=0\},\\
\Zz^\perp {\widetilde V}_p^l(\Tom) &:= \{ u \in  {\widetilde V}_p^l(\Tom): \bl u,v \br_{\omega} =0, \forall v \in  \Zz {\widetilde V}_p^l(\Tom)\}.
\end{align} \end{subequations} 
Then, the discrete Poincar\'e inequalities we want to prove take the following unified form
(where the constant is $\calC=\CP$ or $\CPA$ depending on the context).
\begin{proposition}[Discrete Poincar\'e inequality]\label{th:unifiedDiscretePoincare}
There exists a constant $\calC$, only depending on the mesh shape-regularity parameter $\rho_{\Th}$ and the polynomial degree $p$, such that, for all $\sigma \in \Delta_h$, letting $\omega:= \es(\sigma)$,
\begin{equation} \label{onto_abstract} 
\|u\|_{L^2(\omega)} \le \calC h_\omega \| d^l u \|_{L^2(\omega)},  \qquad \forall u \in \Zz^\perp {\widetilde V}_p^l(\calT_\omega), \quad \forall l\in\{1{:}2\}.
\end{equation}
\end{proposition}
 
\subsection{Reference patches}
For all $\sigma \in \Delta_h$, we 
enumerate the set of vertices in $\Tom$ as $\Vom:=\{ v_0, \ldots, v_{N^{\mathrm v}}\}$,
the set of edges as $\Eom:=\{e_0,\ldots,e_{N^{\mathrm e}}\}$, 
the set of faces as $\Fom:=\{f_0,\ldots,f_{N^{\mathrm f}}\}$, and 
the set of cells as $\Tom:=\{\tau_0,\ldots,\tau_{N^{\mathrm c}}\}$ 
(with $N^{\mathrm c}+1=\rev{\#|\Tom|}$, the \rev{cardinality} of $\Tom$).
Notice that $N^{\mathrm v}$, $N^{\mathrm e}$, $N^{\mathrm f}$, and $N^{\mathrm c}$ 
are bounded from above by a
constant depending only on the mesh shape-regularity parameter $\rho_{\Th}$. 
The topology (and orientation)
of the mesh $\Tom$ is completely described by the connectivity arrays 
\begin{subequations} \begin{align}
&\jev:\{0{:}N^{\mathrm e}\}\times\{0{:}1\}\rightarrow \{0{:}N^{\mathrm v}\}, \\
&\jfv:\{0{:}N^{\mathrm f}\}\times\{0{:}2\}\rightarrow \{0{:}N^{\mathrm v}\}, \\ 
&\jcv:\{0{:}N^{\mathrm c}\}\times\{0{:}3\}\rightarrow \{0{:}N^{\mathrm v}\},
\end{align} \end{subequations} 
such that $\jev(m,n)$ is the global vertex number of the vertex $n$ of the edge $e_m$,
and so on (the local enumeration of vertices is by increasing enumeration order). 
Notice that the connectivity arrays only take integer values and are independent of the
actual coordinates of the vertices in the physical space $\R^3$.

Let $\rho_\sharp>0$ be a positive real number and let $\rev{M_\sharp}$ be a (finite) integer number.
The number of meshes with shape-regularity parameter bounded from above by $\rho_\sharp$ and \rev{cardinality} \rev{bounded from above by} $\rev{M_\sharp}$ with different possible realizations of the connectivity arrays is bounded from above by a constant $\hat N_{\sharp}\eq\hat N(\rho_\sharp,\rev{M_\sharp})$ only depending on $\rho_\sharp$ and $\rev{M_\sharp}$.
Thus, for each $\rho_\sharp$ and $\rev{M_\sharp}$, there is a \emph{finite} set of 
reference meshes, which we denote by $\wTT\eq \wTT(\rho_\sharp,\rev{M_\sharp})$,
such that every mesh $\mathcal{T}$ with the shape-regularity parameter bounded from above by $\rho_\sharp$ and \rev{cardinality} \rev{bounded from above by} by $\rev{M_\sharp}$ has the same connectivity arrays 
as those of one reference mesh in the set $\wTT$. 
We enumerate the reference meshes in $\wTT$ as $\{\wT_1,\ldots,\wT_{\hat N_\sharp}\}$ and fix them once and for all. 
For each reference mesh, the element diameters are of order unity, and the shape-regularity parameter is chosen as small as possible (it is bounded from above by $\rho_\sharp$).
For all $j\in\{1{:}\hat N_{\sharp}\}$, we let $\womega_j$ be the open, bounded, connected, 
polyhedral set covered by the reference mesh $\wT_j$. 

For all $l\in\{1{:}2\}$, we define the piecewise polynomial spaces $V_p^l(\wT_j)$
and $\mV_p^l(\wT_j^{\textsc{a}})$ as in~\eqref{eq:local_spaces}--\eqref{eq:local_spaces_m},
where $\wT_j^{\textsc{a}}$ denotes the Alfeld split of $\wT_j$. 
We set ${\widetilde V}_p^l(\wT_j):=V_p^l(\wT_j)$ or $\mV_p^l(\wT_j^{\textsc{a}})$ depending on the context. 
As above, we use the kernel and orthogonal subspaces such that
\begin{subequations} \begin{align}
  \Zz {\widetilde V}_p^l(\wT_j) &:= \{ \wu \in {\widetilde V}_p^l(\wT_j) : d^l\wu=0\},\\
  \Zz^\perp {\widetilde V}_p^l(\wT_j) &:= \{ \wu \in  {\widetilde V}_p^l(\wT_j): \bl \wu,\wv \br_{\womega_j} =0, \forall \wv \in  \Zz {\widetilde V}_p^l(\wT_j)\}.
\end{align} \end{subequations} 
Using  norm equivalence in finite-dimensional spaces proves the following discrete Poincar\'e inequality in each reference patch and for any $p\ge0$.

\begin{lemma}[Discrete Poincar\'e inequality on reference patches]
For all $j\in\{1{:}\hat N_{\sharp}\}$ and all $p\ge0$, there exists a constant $\CP(j,p)$ such that
\begin{equation} \label{eq:reference_Poincare}
\|\widehat w\|_{L^2(\womega_j)} \le \CP(j,p)  
\| d^l \widehat w \|_{L^2(\womega_j)}, \quad \forall \widehat w \in \Zz^\perp {\widetilde V}_p^l(\wT_j).
\end{equation}
\end{lemma}

\subsection{Piola maps}

Consider an arbitrary $\sigma \in \Delta_h$ with $\omega=\es(\sigma)$ and the corresponding local mesh $\Tom$ such that its shape-regularity parameter is bounded from above by $\rho_\sharp$ and its \rev{cardinality} \rev{is bounded from above} by $\rev{M_\sharp}$. There is an index $j(\Tom) \in \{1{:}\hat N_\sharp\}$ so that $\Tom$ and $\wT_{j(\Tom)}$ share the same connectivity arrays. It follows from~\eqref{eq_rho} that we can take $\rho_\sharp = \rho_{\Th}$. Moreover, as discussed above, the number of tetrahedra in the extended star $\es(\sigma)$ is bounded as a function of $\rho_{\Th}$. Thus, since $\hat N_{\sharp}\eq\hat N(\rho_\sharp,\rev{M_\sharp})$ only depends on $\rho_\sharp$ and $\rev{M_\sharp}$, we infer that $\hat N_{\sharp}$ is bounded \rev{as a} function of $\rho_{\Th}$.

Since $\Tom$ and $\wT_{j(\Tom)}$ share the same connectivity arrays, $\Tom$ can be generated from $\wT_{j(\Tom)}$ by a piecewise-affine geometric mapping
$\bF_{\Tom}\eq\{\bF_\tau:\wtau\rightarrow \tau\}_{\tau\in\Tom}$, where all the geometric mappings $\bF_\tau$ are affine, invertible, with positive Jacobian, and such that $\bigcup_{\tau\in\Tom}\bF_{\tau}^{-1}(\tau)=\wT_{j(\Tom)}$.
For all $\tau\in\Tom$, let $\bJ_\tau$ be the Jacobian matrix of $\bF_\tau$.
We consider the Piola transformations $\psi_{\Tom}^l:L^2(\omega)\rightarrow L^2(\womega_{j(\Tom)})$, for all $l\in\{1{:}3\}$, such that, for all $\tau\in\Tom$,  $\psi_\tau^l\eq\psi_{\Tom}^l|_\tau$ is defined as follows: For all $v\in L^2(\tau)$ or all $\bv\in\bL^2(\tau)$,
\begin{align*}
\psi_\tau^1(\bv) &:= \bJ_\tau^{\textsc{t}} ( \bv \circ \bF_\tau), \\
\psi_\tau^2(\bv) &:= \det(\bJ_\tau) \bJ_\tau^{-1} ( \bv \circ \bF_\tau), \\
\psi_\tau^3(v) &:= \det(\bJ_\tau) ( v \circ \bF_\tau).
\end{align*}
The restricted Piola transformations (we keep the same notation for simplicity)
$\psi_{\Tom}^l : {\widetilde V}_p^l(\Tom) \rightarrow {\widetilde V}_p^l(\wT_{j(\Tom)})$ are isomorphisms. This follows from the fact that $\Tom$ and $\wT_{j(\Tom)}$ have the same connectivity arrays, that $\bF_{\Tom}$ maps any edge (face, tetrahedron) in $\wT_{j(\Tom)}$ to an edge (face, tetrahedron) of $\Tom$ \rev{preserving its orientation}, and that, for each tetrahedron $\tau\in\Tom$, $\psi_\tau^l$ is an isomorphism that preserves appropriate moments \cite[Lemma~9.13 \& Exercise~9.4]{ErnGuermondbook}.
Moreover, the Piola transformations satisfy the following bounds:
\begin{equation} \label{psiBound}
    \|\psi_{\Tom}^l\|_{\calL} \eq \|\psi_{\Tom}^l\|_{\calL(L^2(\omega);L^2(\womega_{j(\Tom)}))} \le C(\rho_{\sharp})(\overline{h}_{\Tom})^{\rev{l-\frac{d}{2}}},
\end{equation}
where $\overline{h}_{\Tom}$ denotes the largest diameter of a cell in $\Tom$. They also satisfy the following commuting properties: 
\begin{equation} \label{eq:commut_Piola}
d^{l}(\psi_{\Tom}^l(v)) = \psi_{\Tom}^{l+1}(d^lv), \qquad \forall v \in \{w\in L^2(\omega):\; 
d^l w\in L^2(\omega)\}.
\end{equation}
We use the shorthand notation $\psi_{\Tom}^{-l}$ for the inverse of the Piola transformations. 
We have 
\begin{equation} \label{psiInverseBound}
    \|\psi_{\Tom}^{-l}\|_{\calL} \eq \|\psi_{\Tom}^{-l}\|_{\calL(L^2(\womega_{j(\Tom)});L^2(\omega))} \le C(\rho_{\sharp})(\underline{h}_{\Tom})^{\rev{-l+\frac{d}{2}}},
\end{equation}
where $\underline{h}_{\Tom}$ denotes the smallest diameter of a cell in $\Tom$. The commuting property~\eqref{eq:commut_Piola} readily gives 
\begin{equation} \label{eq:commut_Piola_inv}
    \psi_{\Tom}^{-(l+1)}(d^l \wv) = d^l(\psi_{\Tom}^{-l}(\wv)), \qquad \forall \wv\in \{\widehat{w} \in L^2(\womega): \; d^l \widehat{w} \in L^2(\womega)\}.
\end{equation}

We will also need the $\R^{3\times 3}$-valued weights 
\begin{equation}
\varrho^1:=\det(\bJ_\tau) (\bJ_\tau^{\textsc{t}}\bJ_\tau)^{-1},
\qquad
\varrho^2:=\det(\bJ_\tau)^{-1}\bJ_\tau^{\textsc{t}}\bJ_\tau = (\varrho^1)^{-1}.
\end{equation} 
Notice that
both weights are piecewise constant on $\wT_{j(\Tom)}$ and take symmetric positive-definite values.
The key property of the weights (e.g. \cite[equation~(18.17)]{ErnGuermondbook}) is that, for all $l\in\{1{:}2\}$,
\begin{equation} \label{eq:Piola_l21}
\bl u,v\br_\omega = \bl \varrho^l \psi_{\Tom}^l(u),\psi_{\Tom}^l(v) \br_{\womega},
\qquad \forall u,v \in \widetilde V_p^l(\Tom).
\end{equation}
We also have the bounds
\begin{equation}\label{eq:varrho_bound}
  \lambda_\flat^l \|\wu\|_{L^2(\womega)}^2 \le \bl \varrho^l \wu,\wu\br_{\womega} \le \lambda_\sharp^l \|\wu\|_{L^2(\womega)}^2,  \qquad \forall \wu\in {\widetilde V}_p^l(\wT_{j(\Tom)}),
\end{equation}
where $0<\lambda_\flat^l\le \lambda_\sharp^l$ denote, respectively, the lowest and largest
eigenvalue of $\varrho^l$ in $\womega$. Invoking, e.g., \cite[equation~(11.3)]{ErnGuermondbook},
we have $\lambda_\flat^1\approx \lambda_\sharp^1 \approx h_\omega$ and
$\lambda_\flat^2\approx \lambda_\sharp^2 \approx h_\omega^{-1}$ (recall that $d=3$), so that
\begin{equation}\label{aux801}
\lambda_\sharp^\ell \lesssim \lambda_\flat^\ell.  
\end{equation}

\subsection{Proof of Proposition~\ref{th:unifiedDiscretePoincare}}

We can now prove the discrete Poincar\'e inequality~\eqref{onto_abstract}.
Let $u \in \Zz^\perp {\widetilde V}_p^l(\Tom)$ 
and let $\widehat q$ be the unique  element in $\Zz^\perp {\widetilde V}_p^l(\wT_{j(\Tom)})$ such that $d^l \widehat q = \psi_{\Tom}^{l+1}( d^l u)$. (The existence and uniqueness of $\widehat q$ follow by reasoning as in~\eqref{eta1def}.)
Using the commuting property~\eqref{eq:commut_Piola}, we infer that
\begin{equation}\label{aux701}
\psi_{\Tom}^{l}(u) - \widehat q \in \Zz {\widetilde V}_p^l(\wT_{j(\Tom)}).   
\end{equation}
The discrete Poincar\'e inequality~\eqref{eq:reference_Poincare} and the fact that  $d^l \widehat q= \psi_{\Tom}^{l+1}( d^l u) $ give
\begin{equation} \label{eq:q_ineq}
\|\widehat q\|_{L^2(\womega)} \leq \CP(j(\Tom),p)\|\psi_{\Tom}^{l+1}\|_{\calL} \| d^lu \|_{L^2(\omega)}.
\end{equation}
Invoking~\eqref{eq:varrho_bound}, we infer that
\begin{equation} \label{eq:bnd1}
\|u\|_{L^2(\omega)}^2 = \| \psi_{\Tom}^{-l} \psi_{\Tom}^{l}(u)\|_{L^2(\omega)}^2 
\leq (\lambda_\flat^{l})^{-1} \|\psi_{\Tom}^{-l}\|_{\calL}^2 \bl \varrho^l\psi_{\Tom}^{l}(u), \psi_{\Tom}^{l}(u)\br_{\womega}.
\end{equation}
Let $\widehat v:=\psi_{\Tom}^{l}(u) - \widehat q$ so that $d^l\widehat v=0$. Invoking the commuting property
for $\psi_{\Tom}^{-l}$ shows that $d^l\psi_{\Tom}^{-l}(\widehat v)=\psi_{\Tom}^{-(l+1)}(d^l\widehat v)=0$,
so that $\psi_{\Tom}^{-l}(\widehat v) \in \Zz {\widetilde V}_p^l(\Tom)$. Since $u \in \Zz^\perp {\widetilde V}_p^l(\Tom)$,
we infer from~\eqref{eq:Piola_l21} that 
\begin{equation}
\bl \varrho^l\psi_{\Tom}^{l}(u), \widehat v\br_{\womega}
= \bl u,\psi_{\Tom}^{-l}(\widehat v)\br_{\omega} = 0.
\end{equation}
Using this identity in~\eqref{eq:bnd1}, we see that
\begin{alignat*}{3}
\|u\|_{L^2(\omega)}^2 
&\le (\lambda_\flat^{l})^{-1}\|\psi_{\Tom}^{-l}\|_{\calL}^2 \bl \varrho^l\psi_{\Tom}^{l}(u), \widehat q\br_{\womega} &&\\
&\leq (\lambda_\flat^{l})^{-1} \lambda_\sharp^l \|\psi_{\Tom}^{-l}\|_{\calL}^2  \|\psi_{\Tom}^{l}\|_{\calL} \|u\|_{L^2(\omega)} \|\widehat q\|_{L^2(\womega)}&& \quad \text{by~\eqref{eq:varrho_bound}}\\
&\leq (\lambda_\flat^{l})^{-1} \lambda_\sharp^l \CP(j(\Tom),p) \|\psi_{\Tom}^{-l}\|_{\calL}^2  \|\psi_{\Tom}^{l}\|_{\calL} \|\psi_{\Tom}^{l+1}\|_{\calL} \|u\|_{L^2(\omega)} \| d^lu \|_{L^2(\omega)} &&\quad\text{by~\eqref{eq:q_ineq}}.
\end{alignat*}
Therefore, using~\eqref{aux801}, \eqref{psiBound}, and~\eqref{psiInverseBound}, we obtain 
\begin{equation}
\|u\|_{L^2(\omega)}\lesssim C(\rho_{\Tom}) \CP(j(\Tom),p) h_{\omega} \| d^lu \|_{L^2(\omega)},
\end{equation}
where we used that $\overline{h}_{\Tom} \leq h_\omega$ and $\overline{h}_{\Tom}/\underline{h}_{\Tom} \le C(\rho_{\Tom})$. We conclude
that~\eqref{onto_abstract} holds true with a constant proportional to
$\max_{j\in\{1{:}\hat N_\sharp\}} \CP(j,p)$, i.e., only depending on $\hat N_\sharp$ and the polynomial degree $p$, thus only depending on the mesh shape-regularity parameter $\rho_{\Th}$ and the polynomial degree $p$.

\bibliographystyle{acm_mod}
\bibliography{references}

\end{document}